\documentclass{article}

\usepackage{arxiv}

\usepackage[utf8]{inputenc} 
\usepackage[T1]{fontenc}    
\usepackage{hyperref}       
\usepackage{url}            
\usepackage{booktabs}       
\usepackage{amsfonts}       
\usepackage{nicefrac}       
\usepackage{microtype}      
\usepackage{lipsum}
\usepackage{amsmath}
\usepackage{inputenc}

\usepackage{subcaption}
\usepackage{rotating}
\usepackage{tikz}
\usepackage{lscape}
\usetikzlibrary{positioning}
\usepackage{array,calc}
\newlength\lena \newlength\lenb \newlength\lenc \newlength\lend
\newcolumntype{P}[1]{>{\centering\arraybackslash}p{#1}} 

\title{Multi-Objective Eco-Routing for Dynamic Control of Connected \& Automated Vehicles}

\author{
	\noindent
  Shadi Djavadian\\
  Laboratory of Innovations in Transportation (LiTrans)\\
  Ryerson University\\
  Toronto, Canada \\
  \texttt{shadi.djavadian@ryerson.ca}
   \And
  Ran Tu\\
  University of Toronto\\
  Toronto, Canada\\
 \texttt{ran.tu@mail.utoronto.ca}
\And   
  Bilal Farooq \\
  Laboratory of Innovations in Transportation (LiTrans)\\
  Ryerson University\\
  Toronto, Canada \\
  \texttt{bilal.farooq@ryerson.ca}
  \And
  Marianne Hatzopoulou\\
  University of Toronto\\
  Toronto, Canada\\  
  \texttt{marianne.hatzopoulou@utoronto.ca} \\
  \\
}

\begin{document}
\maketitle

\begin{abstract}
The advent of intelligent vehicles that can communicate with infrastructure as well as automate the movement provides a range of new options to address key urban traffic issues such as congestion and pollution, without the need for centralized traffic control. Furthermore, the advances in the information, communication, and sensing technologies have provided access to real-time traffic and emission data. Leveraging these advancements, a dynamic multi-objective eco-routing strategy for connected \& automated vehicles (CAVs) is proposed and implemented in a distributed traffic management system. It is applied to the road network of downtown Toronto in an in-house agent-based traffic simulation platform. The performance of the proposed system is compared to various single-objective optimizations. Simulation results show the significance of incorporating real-time emission and traffic state into the dynamic routing, along with considering the expected delays at the downstream intersections. The proposed multi-objective eco-routing has the potential of reducing GHG and NO$_x$ emissions by 43\% and 18.58\%, respectively, while reducing average travel time by 40\%. \\
\vspace{0.5cm}

Reference:\\
Djavadian, S., Tu, R., Farooq, B., Hatzopoulou, M., 2020. Multi-Objective Eco-Routing for Dynamic Control of Connected and Automated Vehicles. Transportation Research Part D: Transport and Environment. 87C: 1-21.
\end{abstract}




\section{Introduction}
\label{S:1}
Transportation systems are one of the major sources of environmental pollution \cite{grote2016including}, accounting for 20\% to 25\% of the world CO$_2$ emission and the world energy consumption \cite{Anagnostopoulou2018}.
Road transportion is also the third leading contributor of NO$_x$, especially from diesel engines. The contribution of local traffic to urban NO$_x$ is estimated to be around 46\% \cite{Transport&Environment}.  
A study by Metrolinx, the public transit authority of the Greater Toronto and Hamilton Area (GTHA), estimated that in 2031 the cost of traffic congestion and environmental pollution to the commuters and economy would reach \$7.8M, and \$7.2M respectively \cite{Metrolinx2008}. Therefore, now more than ever, there is a strong need to mitigate traffic congestion and reduce environmental pollution by taking advantage of new and upcoming technologies.


There are two types of factors that can affect the throughput of traffic network and environmental pollution negatively: a) global factors impacting the entire network (strategic component of driving task, such as route choice, mode choice, and departure time) and b) the local factors that are responsible for local perturbations (manoeuvring and operational levels of driving task) \cite{Michon1985}. The focus of this study is on the global factors in particular route choice of the vehicles. With the advances in communication technologies and automation in vehicles new opportunities have presented themselves to mitigate the negative impact of transportation. In recent years, research on the use of real-time traffic information and dynamic vehicle routing using communication and automation in vehicles has gained interest among researchers as well as traffic engineers. Studies have shown the value of real-time traffic information and dynamic routing in reducing travel time and maximizing capacity \cite{djavadian2018distributed, alfaseeh2018impact}. 

The main objective of the current vehicle routing strategies is predominantly to reduce the travel time, while not considering the byproducts of driving such as environmental pollution. Recently routing Apps like WAZE have adapted a routing scheme that helps drivers avoid certain areas to help reduce emissions in those areas \cite{Sawyer2019}. However, the routing suggestions do not explicitly try to reduce GHG or NO$_x$ emitted per vehicle, meaning the routes are still calculated based on travel time minimization as opposed to emission minimization. Aside from not taking into account the environmental impact, the current vehicle routing systems face additional shortcomings. The current state-of-the-art solutions (e.g. WAZE) are heavily reliant on the information they receive from users as probe vehicles. As such, market penetration rate has a significant impact on the accuracy of route guidance provided. In addition, since the guidance is given by the centralized system, it entails scalability issues. Moreover, same information is provided to all drivers. With drivers' individualistic behaviour and they satisfying their own objective (i.e. fastest route), such routing strategies result in the creation of additional congestion on the network. 

\cite{FarooqBilalandDjavadian} developed a distributed traffic management system based on a network of Intelligent Intersections (I2s) that can dynamically route connected \& automated vehicles using a travel time based routing system called End-to-End Routing for Connected and Automated Vehicles (E2ECAV). Through local level communication among neighbours and frequent dissemination of the local traffic conditions, I2s learn to have the full view of the network conditions. The system needs the vehicles to have Vehicle to Infrastructure (V2I) communication enabled, but automation is not required. \cite{djavadian2018distributed, alfaseeh2018impact, tu2019quantifying} have analyzed the efficiency of the proposed E2ECAV with respect to current solutions on the road network of downtown Toronto in an agent-based traffic simulation and showed that the proposed solution has the potential to reduce travel time as well as environmental pollution, while maximizing throughput. This study adopts the I2s system and extends the single-objective E2ECAV routing to a multi-objective eco-routing where the impacts of both travel time and environmental pollution are explicitly taken into account. To account for the nonlinear speed-emission relationship, highly disaggregate temporal emission models are employed. The following methodological research questions are addressed in this study:

\begin{enumerate}
	\item What criteria are necessary for the comprehensive design of the dynamic distributed multi-objective eco-routing system?
	
	\item How well the proposed multi-objective eco-routing performs when compared to single-objective routing solely based on travel time or emission optimization?
	
\end{enumerate}

The remainder of this paper is organized as follows. In section \ref{S:2} an overview of existing literature on distributed routing and eco-routing is provided. Section \ref{S:3} presents the detailed design of the proposed multi-objective eco-routing strategies. Then a case study of the downtown Toronto network is presented in section 
\ref{S:4} followed by analysis of the results in section \ref{S:5} and concluding discussion in section \ref{S:6}.

\section{Background}
\label{S:2}
To mitigate the negative impact of transportation on the environment and reduce pollution, eco-routing \cite{aziz2012integration} has been considered an alternative to the conventional vehicle routing methods that are based on travel time cost only. Eco-routing is employed to alleviate the unfavourable impact of transportation systems on the environment. Other synonym terms for eco-routing can be: pollution routing \cite{bektacs2011pollution}, or green routing \cite{guo2013evaluation}. The main objective of eco-routing is to route vehicles to the most energy-efficient paths by optimizing their fuel consumption and subsequently reduce environmental emissions (e.g. GHG, NO$_x$). To develop the eco-routing methods, macroscopic and mesoscopic traffic and emission models have been predominantly used. The main justification is the simplicity associated and ease of obtaining the data. However, macroscopic and mesoscopic models are unable to represent the traffic and emissions realistically \cite{alfaseeh2019multi}. For instance, \cite{ guo2013evaluation,patil2016emission,long2016link} depended on average speed to estimate emissions that led to the underestimation. Furthermore, when regression models were employed by \cite{patil2016emission} and \cite{aziz2012integration} to estimate emissions, congested traffic conditions were not taken into account, as regression models are mainly dependent on vehicular speed. Previous studies mostly focused on optimizing a single objective---either GHG or NO$_x$ emission \cite{guo2013evaluation, sun2015stochastic, andersen2013ecotour}. Moreover, in studies that have looked at reducing the negative impact of transportation on the environment, public health is seldom taken into consideration, which is generally represented by NO$_x$ pollutant. In terms of case studies adopted, large scale networks were avoided when the microscopic traffic and emission models were utilized. It is due to the complexity associated and difficulty to obtain microscopic data points. The majority of existing work on the eco-routing employed a centralized control system. However, the advent of intelligent vehicles that can communicate with infrastructure as well as automate the movement provides a range of new options to address key urban traffic issues (e.g. congestion) without a need for a centralized traffic control center. Since early 2000s many researchers have studied the potential for distributed traffic information system e.g. \cite{Kattan2012,Noori2013,AlMallah2017a,AlMallah2017b}. It has been shown that distributed routing overcomes the shortcomings associated with conventional centralized routing systems \cite{FarooqBilalandDjavadian}. These drawbacks are related to the large investment required, high sensitivity to system failure, lack of information relevant to individual trips, and the high complexity when updating the system \cite{yang2006modeling}. Distributed routing systems have recently proved their robustness compared to the centralized solutions \cite{djavadian2018distributed, FarooqBilalandDjavadian}. For a detailed review of existing literature on eco-routing the interested reader is referred to the study conducted by \cite{alfaseeh2019multi}.

\begin{figure}[ht!]
	\centering
	\includegraphics[width=0.75\textwidth]{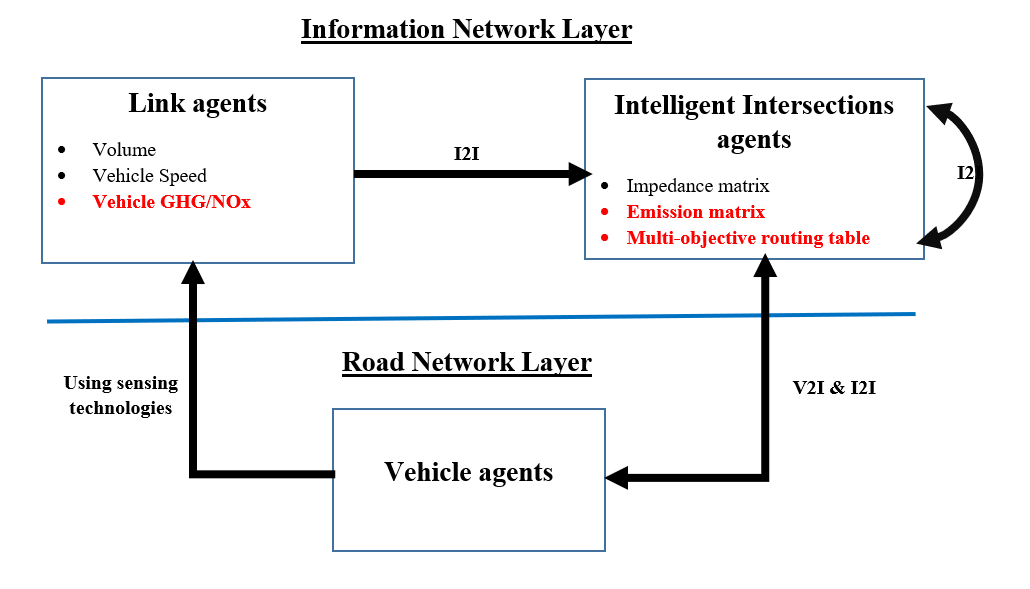}
	\caption{Flow diagram for Multi-Objective E2ECAV routing.}\label{E2ECAV}
\end{figure}

\cite{FarooqBilalandDjavadian} developed a distributed traffic management system design focusing on dynamic End-to-End (E2E) routing based on a network of intelligent intersections (I2) using CAVs. Figure \ref{E2ECAV} presents the schematic view of the I2 system. The framework consists of two network layers and two types of agents. The two network layers are: \emph{Information communication network} where the information exchange between vehicles, links and intersections takes place, and the physical road network where vehicles travel. \emph{Physical road network} is represented by a graph of $G(I, L)$, which consists of I intersections (nodes/vertices) and $L$ links (edges). The two agents are the CAV agents ($v \in V$), and infrastructure agents. There are two types of infrastructure agents, namely: link agents ($l \in L$) and intelligent intersection agents ($I_i^2 \in I^2$). In addition, there are two types of communications in the proposed framework, two-way V2I and I2I taking place over Dedicated Short-Range Communication (DSRC), 5G, or other well-known communication standards.

As shown in Figure \ref{E2ECAV}, real-time traffic information is collected by links using sensors and being frequently exchanged among intelligent intersection using Infrastructure to Infrastructure (I2I) communication. Using the shared information, intelligent intersections very frequently estimate the network traffic state in such a way that they have a consistent, up-to-date, and coherent view, as per the definition of the distributed systems. Each intersection then uses the up-to-date information of the downstream conditions to route the CAVs arriving from upstream toward their announced destination. \textcolor{black}{For further information on the system design, readers are referred to \cite{djavadian2018distributed,DjavadianS.FarooqB.VasquezR.Yip2019, alfaseeh2018impact, tu2019quantifying}.}

\cite{djavadian2018distributed} showed that the proposed E2ECAV is capable of reducing travel time while maximizing network throughput, especially under high market penetration rate (MPR) and high congestion level \cite{alfaseeh2018impact}. \cite{tu2019quantifying} assessed the impact of E2ECAV (with minimizing only the travel time (TT)) on traffic emissions (GHG and NO$_x$) under different MPR  of CAVs and different traffic congestion levels. It was found that higher MPRs of CAVs contributed to significant GHG emission reduction on the network in comparison to \textcolor{black}{ HDV pre-trip dynamic routing with no en-route routing}. However, as NO$_x$ is more sensitive to aggressive driving of higher speed range, reduction in NO$_x$ happens for specific MPRs of CAVs under certain traffic conditions.
Although the results from \cite{tu2019quantifying} showed that E2ECAV with a single-objective (i.e. travel time) optimization is capable of reducing air pollution, the effect of traffic emissions such as GHG and NO$_x$ were not explicitly included in the objective function. We hypothesize that explicitly accounting for such emissions in the objective function will lead to even further reductions. \cite{alfaseeh2019eco} showed that by including such factors in the objective function can result in the improvement in total emissions from network. However, the study used a myopic objective function and did not account for the downstream link characteristics i.e. travel time and delay. 

\section{Methodology}
\label{S:3}
This study adopts the I2 system \cite{FarooqBilalandDjavadian} and extends the routing to multi-objective eco-routing (refer to red features added in Figure \ref{E2ECAV}), where the impact of both travel time and environmental pollution are taken into account explicitly. In real life scenario for multi-objective eco-routing, link agents collect real-time emission data using related sensors such as portable emission measurement system (PEMS). In this case study we use MOtor Vehicle Emission Simulator (MOVES) to emulate those emission senors. MOVES determines an operating mode (opmode) for each second based on the instantaneous speed and acceleration. It then assigns an emission rate to the opmode  by vehicle characteristics \cite{tu2019quantifying, Liu2019, UnitedStatesEnvironmentalProtectionAgency2019}. A vehicle type (passenger car or truck powered by gasoline) with a model year (1988-2018) is randomly assigned to each vehicle in the simulation based on the on-road vehicle fleet composition \cite{MTO2016}. \textcolor{black}{In this study, the input data for estimating emissions are obtained from a microscopic simulation. The MOVES model is calibrated based on the emission measurement from multiple sources, including dynamometer (in laboratory) and PEMS (on-road). }

The agent-based microsimulation used in this study is adapted from \cite{djavadian2018distributed}, which uses the intelligent driver car following model \cite{Treiber2000} and has a time-step of 1sec for the generation of instantaneous vehicle operating status. In the proposed system (Figure \ref{E2ECAV}), link infrastructure collects real-time traffic and emission data and frequently exchange them with the downstream intersections. The report sent to each downstream intersection may contain the following information: the space mean speed $U_{(l,{\Delta}_{j})}$ on link $l$ at time interval ${\Delta}_{j}$ and ${{{\hat{\text{CO}}^{+}}_{{2}_{(l,{\Delta}_{j})}}}}$ which is space mean GHG (in ${\text{CO}}_{2-eq}$) emission rate ($gram/veh/s$) of link $l$ at time interval ${\Delta}_{j}$ (calculated by Equation \ref{spacemeanemission}).

\begin{equation}
{{{\hat{\text{CO}}_{{2}_{(l,\Delta_j)}}}^+} =\displaystyle \frac{ \displaystyle \sum_{p=1}^{P} \frac{ \displaystyle\sum_{\omega=1}^{\Omega} {{\hat{\text{CO}}_{{2}_{(p,l,\Delta_\omega)}}}^+}}{\Omega}}{P}}
\label{spacemeanemission}
\end{equation}

$\displaystyle {{\hat{\text{CO}}_{{2}_{(p,l,\Delta_{\omega})}}}^+}$: space mean emission \textcolor{black}{($gram/veh/s$)} per vehicle for section ($p$) during interval $\Delta_\omega$; each link $l$ is divided into ($p$) sections \textcolor{black}{in order to capture different traffic regimes on the link.} It is calculated using Equation \ref{spacemeanemission2}.

\begin{equation}
\displaystyle{\hat{{\text{CO}}^+_{2_{(p,l,\Delta_\omega)}}}} = \frac{\displaystyle \sum_{m=1}^{M} {{\text{CO}^+}_{{2_{(m,p,l,\Delta_{\omega})}}}}}{M}
\label{spacemeanemission2}
\end{equation}

M: is total number of vehicles on segment $p$ of link $l$ during time interval $\Delta_\omega$, $M \subset V$.

${{\text{CO}_2}+}_{(m,p,l,\Delta_\omega)}$: is the emission of vehicle $m$ ($m \in M$) on segment $p$ of a link $l$ during time interval $\Delta_\omega$ calculated using second-by-second acceleration and speed data along with MOVES. In addition to GHG (in $\text{CO}_2$-eq), MOVES also estimates NO$_x$.

Using the updated real-time traffic and emission information received from adjacent link agents, I2s guide CAVs through the network based on a predefined routing objective, e.g. minimizing individual travel time. The general objective function for the multi-objective routing is presented by Equation \ref{e:1}.

\begin{eqnarray}
\min  \sum_{l} \alpha_l\times[\beta_T\times{\hat{T}_{(l,\Delta_j)}}+ \beta_{\Pi}\times{\hat{\Pi}_{(l,\Delta_j)} } + \beta_{\text{CO}_2}\times{{{\hat{E}}_{\text{CO}_{{2}_{(l,\Delta_j)}}}}}  ]
\label{e:1}
\end{eqnarray}

${{\hat{T}}_{(l,\Delta_{j})}}$ is the average travel time on link $l$ at time interval $\Delta_j$. ${\hat{E}_{\text{CO}_{{2}_{(l,\Delta_j)}}}}$ is the average GHG emission (gm) per vehicle on link $l$ at time interval $\Delta_j$.   $\hat{\Pi}_{(l,\Delta_j)}$ is the average \textcolor{black}{idling penalty} at the downstream intersection of link $l$ at time interval $\Delta_j$. \textcolor{black}{For the purpose of this study the idling penalty is calculated in terms of the average time a vehicle spends right at the intersection before it obtains right of way to cross the intersection, which is affected by the traffic volume on other approaches.} The \textcolor{black}{decision} variable is $\alpha_l$, which is a binary variable, if link $l$ is selected $\alpha_l=1$ otherwise is 0. $\beta_{TT}, \beta_\Pi$, and $\beta_{\text{CO}_2}$ are binary routing strategy parameters, if a certain cost is considered in the objective function $\beta=1 $ otherwise 0. Two different eco-routing strategies are proposed and the $\beta$ values are adjusted accordingly. The significant difference between the two proposed strategies is that one only aims to minimize GHG emission (single objective) on the network where as the other takes into consideration other factors such as travel time and \textcolor{black}{idling penalty} in addition to GHG emission (multi-objective). The aim is to show the benefit of multi-objective eco-routing strategies and the inter-connectivity between GHG emission, \textcolor{black}{idling penalty} and travel time. It is worth mentioning that although NO$_x$ is not directly included in Equation \ref{e:1}. By minimizing travel time and at the same time GHG emission, it indirectly finds a reasonable compromise in terms of reducing NO$_x$ emission.   

In the context of representing emissions in the objective function, two routing strategies are used:

\begin{itemize}
	
	\item [R1.]  \textbf{$GHG^+_{raw}$} ($\beta_{\text{CO}_2}=1$): minimizing the individual vehicle's GHG emission where link cost is based on raw GHG emission, taking into account the individual link's ($l$) physical (e.g. length) and traffic characteristics (e.g. average speed).
	
	\item[R2.]  \textbf{$TT^* + GHG^+_{raw}$} ($\beta_{\Pi}=1$, $\beta_{\text{CO}_2}=1$, and $\beta_{TT}=1$): multi-objecitve eco-routing based on link GHG emission,  average travel time and average \textcolor{black}{idling penalty}. 
	
\end{itemize}

\textcolor{black}{It should be mentioned that the two strategies defined in this paper are agnostic and can be implemented on any other routing algorithms/schemes aside from the E2ECAV routing framework.}

\subsection{$\text{GHG}^+_{\text{raw}}$ Routing Strategy (R1)}
The routing strategy presented in this section ( ${GHG}^+_{\text{raw}}$) in addition to second-by-second vehicle emission estimated by MOVES, also takes into account both physical and real-time traffic state characteristics of each link when calculating link costs. In addition, the link GHG is calculated based on space mean GHG taking into account spatio-temporal distribution of the vehicles and different traffic regimes on the link. The objective function for $GHG^+_{raw}$ routing is presented by Equation \ref{e:R1}, which is a reduced version of Equation \ref{e:1}, where $\beta_{\text{CO}_2}=1$ and other $\beta$ values are set to 0. ${\hat{E}_{\text{CO}_{{2}_{(l,\Delta_j)}}}}$ is average emission ($gm$) of $\text{CO}_2$ per vehicle on link ($l$) at time interval $\Delta_j$ calculated as shown by Equation \ref{e:E}. $\hat{T}_{(l,\Delta_j)}$ is average travel time ($min$) of link $l$ at time interval $\Delta_j$ and is calculated using Equation \ref{averagetraveltime_eq}. $\Phi_{l}$ is length of link $l$ ($km$). \textcolor{black}{It should be noted that $\hat{T}_{(l,\Delta_j)}$ implicitly takes into account the wait time on link $l$ at time interval $\Delta_j$.}

\begin{equation}
\displaystyle
R1= \min \sum^{l} \Big( \alpha_{l}\times{\hat{E}_{\text{CO}_{{2}_{(l,\Delta_{j})}}}} \Big)
\label{e:R1}
\end{equation}

\begin{equation}
\displaystyle
\hat{E}_{{\text{CO}}_{2_{(l,\Delta_{j})}}}={\hat{{\text{CO}}}_{2_{(l,\Delta_{j})}}^{+}}\times \hat{T}_{(l,\Delta_{j})}
\label{e:E}
\end{equation}

\begin{equation}
\hat{T}_{(l,\Delta_j)}=\frac{\Phi_{l}} {U_{(l,\Delta_j)}}
\label{averagetraveltime_eq}
\end{equation}

\subsection{$\text{TT}^* + \text{GHG}^+_{\text{raw}}$ Routing Strategy (R2)} 
\textbf {R1} implicitly takes into account the cost of travel time on link ($l$) at time interval $\Delta_j$. However, under the routing strategy \textbf{R2} ($TT^* + GHG^+_{raw}$), both the average link travel time (i.e. $\alpha\times\hat{T}_{(l,\Delta_j)}$) and the average \textcolor{black}{idling penalty (min)} (i.e. $\hat{\Pi}_{(l,\Delta_j)}$) at downstream intersection are explicitly taken into account.  
The multi-objective eco-routing is presented by Equation \ref{e:R2}. 

\begin{equation}
\begin{split}
	R2&= \min  \sum_{l} \alpha_{l}\times[W_{T}\times{\hat{T}_{(l,\Delta_{j})}}+\\ &W_{\Pi}\times{\hat{\Pi}_{(l,\Delta_{j})} } +\\ &W_{\text{CO}_{2}}\times{\hat{E}_{\text{CO}_{2_{(l,\Delta_{j})}}}} 
\end{split}
\label{e:R2}
\end{equation}

$W_{T}$, $W_{\Pi}$ and $W_{\text{CO}_2}$, are the weights used when a multi-objective routing scenario is applied. For the purpose of this study, all the objectives are converted to minutes, as such $W_{T}$, $W_{\Pi}$ are set to one, whereas $W_{\text{CO}_2}$ is the monetary value of GHG converted to minutes. The monetary value of GHG emission is \$15.77/Ton based on Qu\'ebec's Cap-and-Trade program \cite{GHG_cost}. The \$0.35/min is used to convert the monetary value of GHG to the equivalent minute units \cite{usdot2014}.


For simplicity, from this point forward, the travel time objective function without \textcolor{black}{idling penalty} (Equation \ref{e:TT}) will be referred to as $TT$ and travel time objective function with the \textcolor{black}{idling penalty} (Equation \ref{e:R2reduced}) will be referred to as $TT^*$.

\begin{eqnarray}
TT=\min  \sum_{l}^{} {\alpha\times\hat{T}_{(l,\Delta_j)}}
\label{e:TT}
\end{eqnarray}

\begin{eqnarray}
TT^*= \min  \sum_{l}^{} {\alpha\times(\hat{T}_{(l,\Delta_j)} + \hat{\Pi}_{(l,\Delta_j)})}
\label{e:R2reduced}
\end{eqnarray}

The aim of the routing strategy R2 is to route the vehicles in such a way that not only their travel time is minimized, but also their GHG emission and the \textcolor{black}{idling penalty} at the intersections. Routing strategy R2, takes into consideration both user and system welfare. 

\begin{figure}[!ht]
	\begin{subfigure}{0.5\textwidth}
		\includegraphics[width=\linewidth, height=4.5cm]{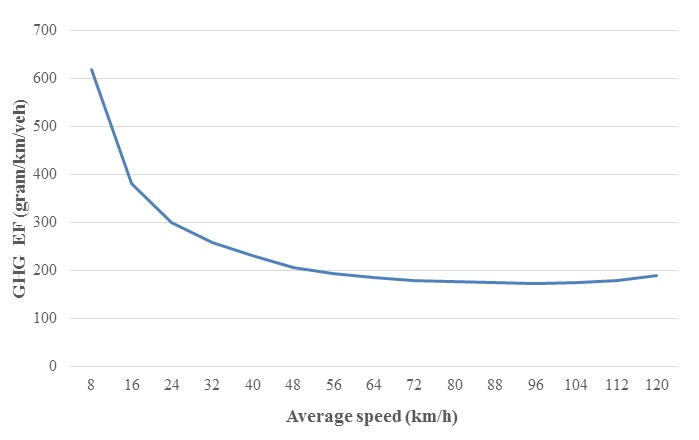} 
		\caption{\textcolor{black}{GHG emission factor}}
		\label{GHG EF}
	\end{subfigure}%
	\begin{subfigure}{0.5\textwidth}
		\includegraphics[width=\linewidth, height=4.5cm]{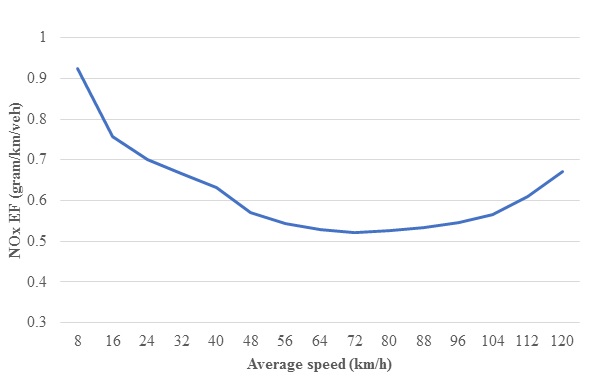}
		\caption{\textcolor{black}{NO$_x$ emission factor}}
		\label{NO$_x$ EF}
	\end{subfigure}%
	\caption{\textcolor{black}{Relationship between GHG and NO$_x$ emission factors and average vehicle speed}}
	\label{GHGEF and NOXEF}
\end{figure}

\textcolor{black}{Although NO$_x$ is not directly included in the objective function of R2, its emission level is still indirectly affected. As shown in Figure \ref{GHGEF and NOXEF} the relation between GHG and NO$_x$ and speed are quasi-convex and non-monotonic. From Figure \ref{GHGEF and NOXEF} it can be seen that both GHG and NO$_x$ are sensitive to low speed, meaning that the longer the vehicle spends on the network, it contributes more to NO$_x$ and GHG. As such there is a need for congestion mitigation and Equation \ref{e:R2}, by minimizing the travel time, ensures that vehicles don't stay on the network for too long. Although both GHG and NO$_x$ are sensitive to high speed levels, as can be seen from Figure \ref{NO$_x$ EF}, the susceptibility of NO$_x$ is higher. Based on the information presented in Figure \ref{GHGEF and NOXEF} that the optimal emission (${E}^{*}$) is not at free flow speed. Moreover, it is drastically affected by acceleration and deceleration on a link. By minimizing GHG, strategy R2 ensures that vehicles don't travel at very high speed. Hence reducing GHG levels and indirectly finding a reasonable compromise in terms of reducing NO$_x$ emission. In general, when GHG is optimized, it ensures that the speeds on the links are not too high and not too low. It is worth mentioning that data for creating Figure \ref{GHGEF and NOXEF} is obtained from MOVES and that Figure \ref{GHGEF and NOXEF} is for illustration purpose.}

\section{Case Study}
\label{S:4}

We used the road network and existing demand from downtown Toronto, Canada. One of the main reasons for selecting this case study was to test the proposed algorithm on a network with high current demand. This network faces recurrent congestion during the morning and afternoon peak periods. It is partly due to the high vehicle ownership and partly due to the under-utilization of some roads and over-utilization of others. It was assumed that the vehicles could operate independently before entering the downtown core, but once they are in the downtown core, they are required to yield the control to the E2ECAV system.

\begin{figure}[!h]
	\begin{center}
		\includegraphics[width=6.5in]{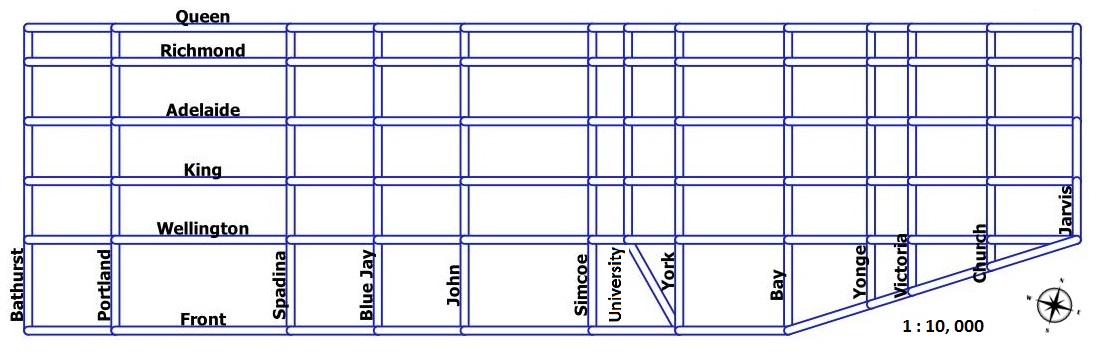}
		\caption{\textcolor{black}{Downtown Toronto street network.}}
		\label{simplifiedtorontonetwork}
	\end{center}
\end{figure}

Figure \ref{simplifiedtorontonetwork} presents the test network (700m x 2,615m), which consists of 76 nodes and 223 links. \textcolor{black}{The network characteristics are presented in Table \ref{network}. To be able to compare, we used the same network as in \cite{tu2019quantifying}. However, it is worth mentioning that E2ECAV routing framework is transferable to larger networks, as shown by \cite{Meshkani2019}, where the framework was applied to an expanded downtown Toronto network with 268 nodes and 839 links.}   In this study, we simulated the current demand during the morning peak period arising from vehicles whose origins and/or destinations are in the downtown core. The time-dependent zonal demand with 5min interval is obtained by applying 2018 growth factor on the travel data from the 2011 Transportation Tomorrow Survey (TTS) of Toronto. The demand is randomly loaded on the nearest intersections. During the study period (7:45am-8:00am) 3,477 vehicles used the network. Within 5min, the demand is then randomly distributed using Poisson arrivals. The simulation time step is 1sec, while for the information communication network, an update interval ($\Delta_j$) of 1min and an intermediate update interval ($\Delta_\omega$) of 20 seconds are used. The simulation is stopped when the last vehicle reaches its destination---thus creating a common base for comparison. A 5min warm-up period is also used. The movements of the vehicles at each intersection are modelled using first in first out (FIFO) strategy, but any other scheme can also be implemented. 

\begin{table}[h!]
	\caption{\textcolor{black}{Network characteristics}}
	\label{network}
	\centering
	\small
	\begin{tabular}{|c|c|c|c|}
		\cline{1-4}
		\textbf{Street Name}                     & \textbf{Speed} & \textbf{Num. Lanes} & \textbf{Direction} \\
		\hline
		Bathurst                         & 60    & 2          & two-way   \\
		\hline
		Portland                         & 40    & 1          & two-way   \\
		\hline
		Spadina                          & 60    & 4          & two-way   \\
		\hline
		Blue Jays                        & 40    & 2          & two-way   \\
		\hline
		John                             & 40    & 2          & two-way   \\
		\hline
		Simcoe                           & 40    & 2          & one-way (SB)   \\
		\hline
		University                       & 60    & 3          & two-way   \\
		\hline
		York                             & 40    & 4          & one-way (NB)   \\
		\hline
		Bay                              & 60    & 2          & two-way   \\
		\hline
		Yonge                            & 60    & 2          & two-way   \\
		\hline
		Victoria                         & 40    & 2          & two-way   \\
		\hline
		Victoria between King \& Front    & 10     & 1          & one-way (NB)   \\ 
		\hline
		Church                           & 40    & 2          & two-way   \\
		\hline
		Jarvis                           & 60    & 2          & two-way   \\
		\hline
		Front                            & 80    & 2          & two-way   \\
		\hline
		Front between York \& University & 40    & 1          & two-way   \\
		\hline
		Wellington                       & 40    & 2          & one-way (EB)  \\
		\hline
		King                             & 60    & 2          & two-way   \\
		\hline
		Adelaide                         & 60    & 4          & one-way (EB)   \\
		\hline
		Richmond                         & 60    & 3          & one-way (WB)  \\
		\hline
		Queen                            & 60    & 2          & two-way  \\
		\hline
	\end{tabular}
\end{table}


The main objective here is to evaluate the design and performance of the distributed multi-objective eco-routing using a range of scenarios. The performance measures used are: average travel time, Vehicle Kilometre Travelled (VKT), network-level GHG and NO$_x$ emissions. In total, five scenarios are tested as presented in Table \ref{table1} to show the significance of incorporating travel time, \textcolor{black}{idling penalty} and GHG emission in the routing objective function.  For all scenarios except 1, it is assumed that all the vehicle are CAVs. In scenario 1, all vehicles are Human Driven Vehicles (HDVs) and are routed based on pre-trip dynamic shortest path at the time of their entry into the network. \textcolor{black}{As shown in \ref{hdv_routing}, we compared various routing strategies for HDV and found that pre-trip dynamic routing performed the best for HDV. Under pre-trip dynamic routing, pre-trip shortest path is selected at the time of entry in the network based on current traffic conditions with no en-route routing. The current traffic condition is estimated by a central controller having full view of the network and access to frequently updated information.}

\begin{table}[!h]
	\caption{Routing scenarios investigated}
	\label{table1}
	\begin{center}
		\begin{tabular}[!t]{l c c c c}
			\hline
			Scenario\# & Routing & Routing Type & obj. Eq'n \\
			& obj.\\
			\hline
			\hline
			
			$S1$ & TT &  Pre-trip Dynamic & \ref{e:TT}   \\
			$S2$ & TT  & E2ECAV & \ref{e:TT}   \\
			$S3$ & R1: GHG  & E2ECAV & \ref{e:R1}    \\
			$S4$ & TT* & E2ECAV & \ref{e:R2reduced}   \\
			$S5$ & R2: TT* + GHG & E2ECAV & \ref{e:R2} \\

			\hline
		\end{tabular}
	\end{center}
	\label{scenarios investigated}
\end{table}

\textcolor{black}{
	\subsection{Implementation of vehicle type}
	As discussed in previous section, two types of vehicles are considered in this case study, HDVs and CAVs. The HDVs and CAVs differ from each other in the simulation in three different ways:
	\begin{enumerate}
		\item [1.] \textbf{Routing}: HDVs are provided with pre-trip dynamic route at the time of entry in the network. They follow the same path to their destination and no en-route routing is available. Whereas, CAV agents solely rely on the instructions given to them by the I2s for routing. Once on a particular link $l$, CAV agents via V2I communicate their location and speed to the link agent ($l$). Similarly, once a CAV agent arrives at an intersection $I_i^2$, it communicates its destination to the $I_i^2$, who routes the CAV to the next $I_{(i+1)}^2$ on the path that meets the objective of the controller ($I_i^2$). 
		\item[2.] \textbf{Communication}: HDVs do not have any active communication ability and they are responsible for all driving tasks. Whereas, CAVs are intelligent vehicles that are capable of sensing, actuating and communication. CAVs are responsible for all movement related tasks except routing.
		\item[3.] \textbf{IDM \cite{Treiber2000} car following parameters}: The reaction time and minimum safe distance parameters for the CAVs are set to be half of the HDVs. The rationale behind this is that these intelligent vehicles can perceive the changes in the environment quicker than humans and can maintain the acceleration and deceleration more smoothly. Maximum acceleration and deceleration is the function of the engine technology and vehicle weight, so in our study, they are the same for both vehicle types.
		\\
\end{enumerate}}



\section{Results and Discussion}
\label{S:5}
Here we present a detailed analysis of the simulation results for scenarios listed in Table \ref{table1}.

\subsection{Network level travel time, VKT, GHG, and NO$_x$}
\label{Travel time, VKT, GHG, and NO$_x$}
\begin{figure}[!ht]
	\begin{subfigure}{0.5\textwidth}
		\includegraphics[width=\linewidth, height=4.5cm]{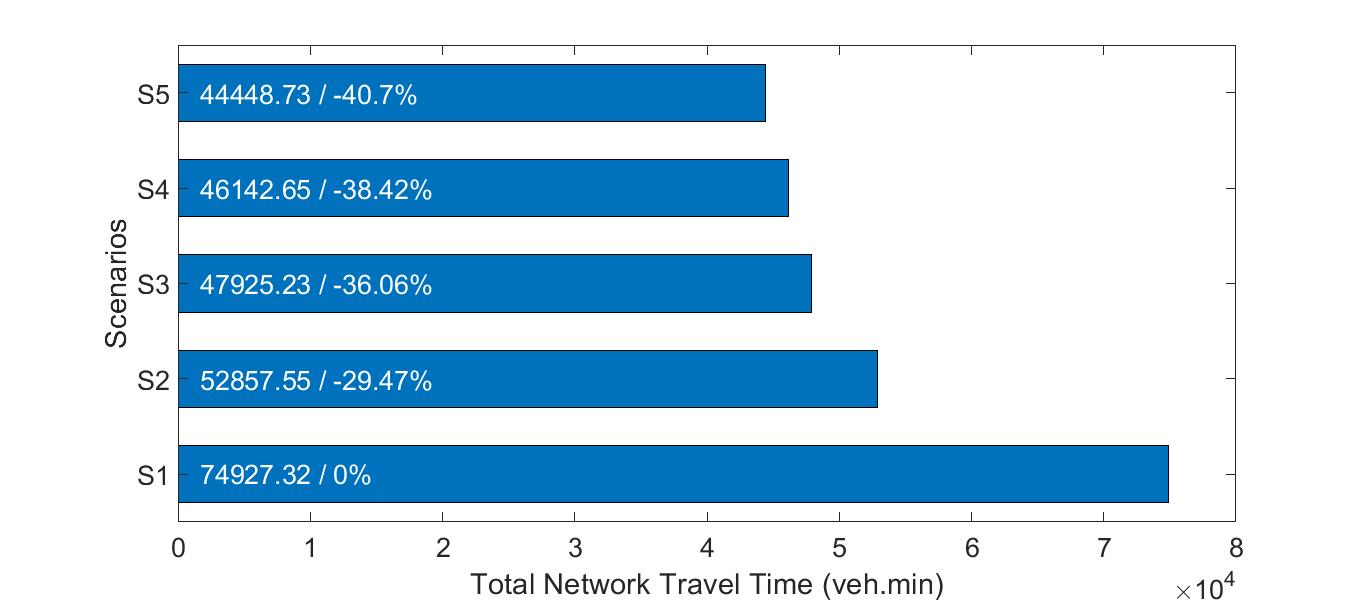} 
		\caption{\textcolor{black}{Total network} TT}
		\label{TT}
	\end{subfigure}%
	\begin{subfigure}{0.5\textwidth}
		\includegraphics[width=\linewidth, height=4.5cm]{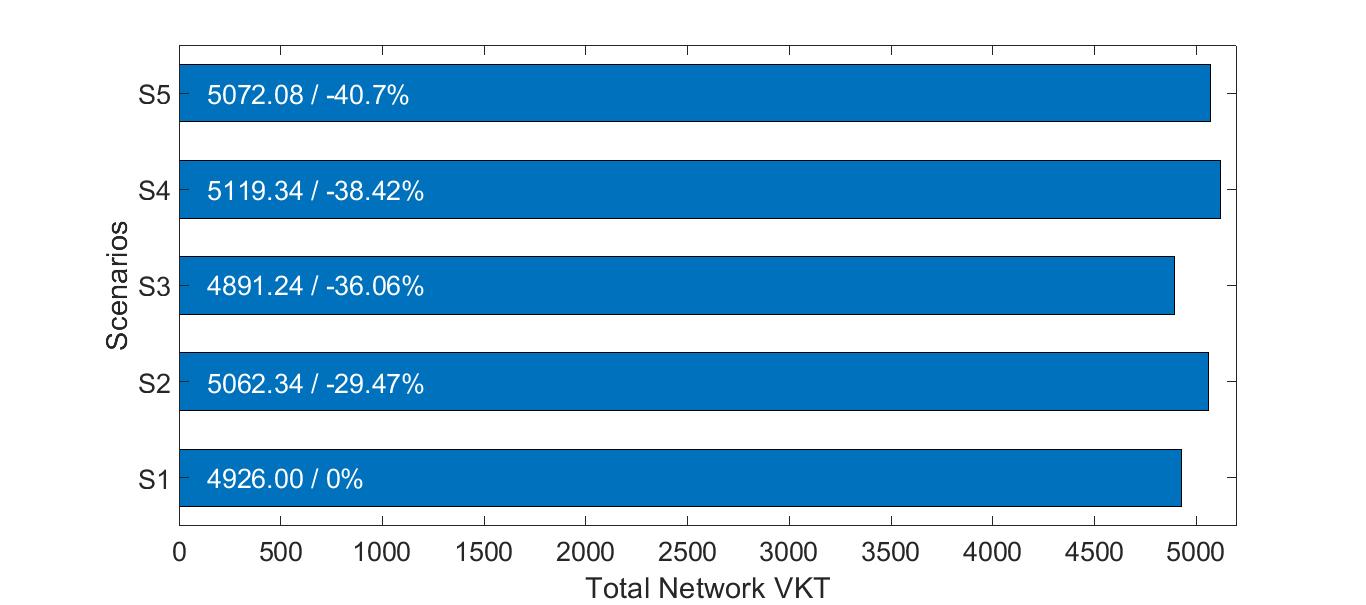}
		\caption{\textcolor{black}{Total network} VKT}
		\label{VKT}
	\end{subfigure}%
	\caption{\textcolor{black}{Total network} TT and VKT of the investigated scenarios}
	\label{TT and VKT}
\end{figure}

From Figure \ref{TT}, it can be observed that scenario S5 resulted in 40.7\% decrease in travel time when compared to scenario S1. A major reason for this significant improvement is that S5 employs E2ECAV, a distributed routing system, where I2s are constantly aware of the real-time traffic condition on the network and are able to periodically re-distribute traffic to the less congested parts of the network, resulting in lower travel times. As such, all scenarios that employed E2ECAV (S2--S5) resulted in a significant reduction in travel time when compared to S1. 
Under S5, the average \textcolor{black}{idling penalty} at the downstream intersection of each link $l$ was taken into consideration (refer to Equation \ref{e:R2}), which resulted in vehicles avoiding links with longer queues. Hence, preventing longer delays and dissipating existing queues faster. Routing under S4 with $TT^*$ as its objective function (Equation \ref{e:R2reduced}) was able to outperform S1--S3 by 38.42\%, 9\%, and 2\%, respectively. Aside from taking into account the current link travel time and average \textcolor{black}{idling penalty} information, the objective function employed in S5 also takes into consideration the minimization of GHG emissions. Thus, resulting in further travel time reduction because having GHG as an objective ensures the balance in distributing traffic on the network and better utilizing the available infrastructure by maintaining optimal speed on each link (not too high and not too low). The aforementioned characteristic of GHG routing, is the reason behind S3's ability to outperform S1 and S2 as well as comparatively reduce network level travel time by 36 \% and 7\%, respectively. \textcolor{black}{ The difference between scenario S2 and multi-objective routing scenario S5 is analogous to that of user equilibrium and system optimal. In user equilibrium the sole objective is to minimize users' travel time (TT). Whereas, in system optimal the objective is to minimize total system travel time by minimizing marginal cost which consists of user costs and additional social cost. In the multi-objective routing scenarios S4 and S5 aside from travel time, we also included idling penalty and GHG emission that are social cost of driving. If all vehicles are sent to shortest path (S2) it leads to congestion formation on shortest paths, however S4 and S5 by taking into consideration the idling penalty and GHG emission put in place restriction on the number of vehicles sent to shortest paths.}

The results from VKT analysis are presented in Figure \ref{VKT}. Overall the mean VKT for the five scenarios appears to be the same (1.4km). However, even-though en-route routing strategies based on E2ECAV implementation (S2-S5) were able to reduce travel time significantly compared to dynamic pre-trip routing under HDV scenario (S1), except in the case of S3 this decrease in travel time came at the expense of slight increase in VKT, approximately $ 3\%$ (max). The reason for this minuscule increase in VKT under E2ECAV is due to CAVs being re-routed en-route and distributed over the network to avoid congestion. In the case of routing scenarios that had travel time explicitly as their objective, increase in VKT did not necessary mean an increase in the travel time. This could be because a longer route may have higher speed, thus a shorter travel time. On the other hand, in the case of GHG emission, longer trips meant higher GHG emission. Therefore, the objective function was more mindful of VKT increase. As can be seen in S3, where the routing is based on solely minimizing GHG emissions, the VKT is lowest when compared to the other four scenarios. \par
\begin{figure}[!ht]
	\begin{subfigure}{0.5\textwidth}
		\includegraphics[width=\linewidth, height=4.5cm]{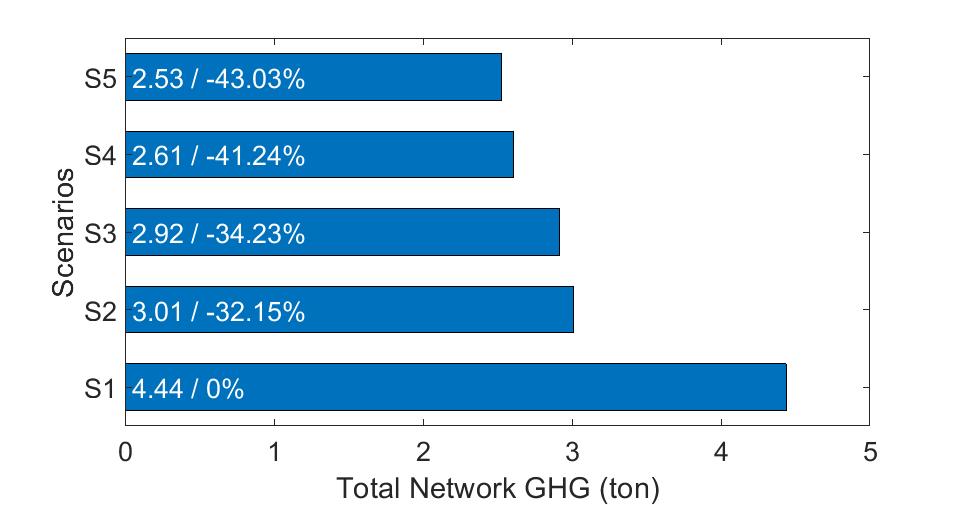} 
		\caption{Total Network GHG }
		\label{GHG}
	\end{subfigure}%
	\begin{subfigure}{0.5\textwidth}
		\includegraphics[width=\linewidth, height=4.5cm]{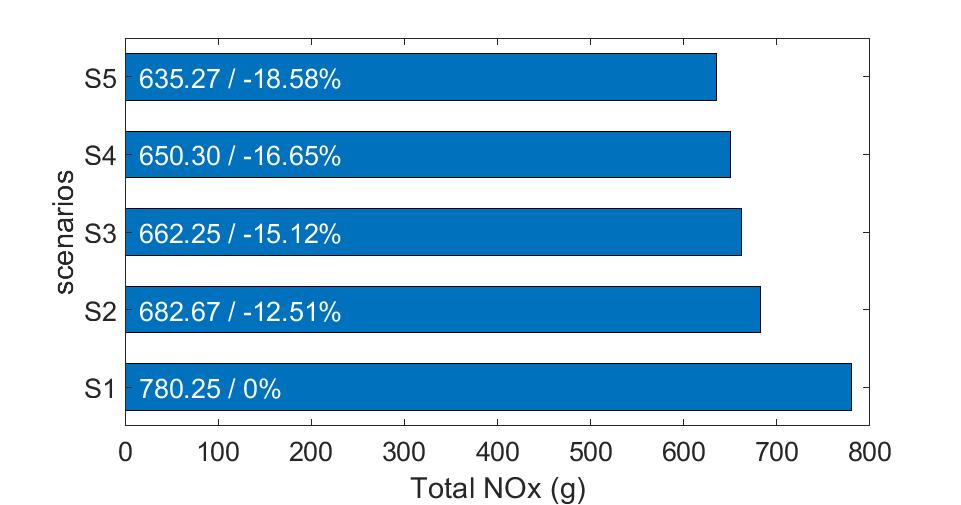}
		\caption{Total Network NOx}
		\label{NOx}
	\end{subfigure}%
	\caption{Total GHG and NO$_x$ of the investigated scenarios}
	\label{GHG and NOx}
\end{figure}

The comparative analysis of GHG emissions for the five scenarios is presented in Figure \ref{GHG}. The reduction in GHG emission levels follows the same pattern as the travel time reduction with scenario S5 having the most significant decrease in GHG emission levels (43\%) when compared to S1. The reason behind this noticeable improvement in GHG emission under S5 is that communicating intelligent intersections are constantly aware of both the real-time traffic condition and GHG emission levels on the network. Based on which, they are able to periodically re-distribute traffic to parts of network with the aim to not only decrease travel time, but at the same time decrease GHG emission levels. There are three main factors that can affect the GHG emission level on the network:
\begin{itemize}
	\item[a)] The amount of time vehicles spend on the network has a significant effect on the emission level. The less time they spend on the network, the less emission they produce. This is why scenarios S2--S5 that are based on more efficient routing scheme, were able to manage GHG emission levels better than S1.
	\item[b)] Although shorter travel time means lower GHG emissions level, just reducing travel time by itself does not guarantee optimal GHG emission level. Shorter travel time may also mean higher speed. Due to the quasi-convex relationship between GHG and speed, higher speed does not necessarily mean lower GHG emission. To reduce the GHG emissions, it is necessary to maintain network speed at optimal emission speed (not too high and not too low). In Figure \ref{GHG}, it can be observed that scenario S3 has the ability to reduce GHG emission more than S2. This reduction can be associated to the fact that the objective function used in S3 (Equation \ref{e:R1} ) takes into consideration both GHG and travel time reductions simultaneously, albeit travel time reduction being secondary to GHG reduction and only being taken into consideration implicitly.
	\item[c)] As stated by Natural Resources Canada, a simple and effective way to manage and reduce GHG emission levels and reduce the production of $\text{CO}_2$ emissions from light-duty vehicles is to eliminate vehicle idling. It is estimated that if motorists can avoid idling by as little as three minutes every day, collectively, they can reduce $\text{CO}_2$ emissions by as much as 1.4 million tonnes annually \cite{Canada2016}. From the results obtained it can be seen that S4 was able to reduce GHG emission by 41\% compared to S1  by 9\% and 7\% compared to S2. This noticeable reduction in GHG emissions occurred due to the fact that as illustrated by Equations \ref{e:R2reduced}, S4 in addition to link travel time, explicitly takes into account average \textcolor{black}{idling penalty} at the downstream intersection of link.
\end{itemize}

Scenario S5 with Equation \ref{e:R2} as its objective function has the ability to outperform S1--S4 because it addresses the three aforementioned factors that affect GHG emissions level. It simultaneously and explicitly, takes into account travel time, average \textcolor{black}{idling penalty} and GHG emission level. R2 aims to route vehicles in such a way that it not only minimizes travel time, but also ensures minimum idling and minimum GHG emission by maintaining optimal emission speed. 

Figure \ref{NOx} presents the results of NO$_x$ emissions level analysis for the five scenarios. It can be seen that although none of the scenarios S2--S5 explicitly considered NO$_x$ in their respective objective functions, they all resulted in a significant decrease in NO$_x$ emissions level when compared to S1, with S5 achieving the most reduction of 18.58\%. Scenarios S2 and S4 were able to achieve reduction in NO$_x$ level by ensuring that vehicles spend less time on the network and minimize the idling time. Whereas, S3 and S5 were able to achieve NO$_x$ reduction by taking into account other environmental byproducts of driving such as GHG emissions.

\subsection{GHG, NOx and Speed over time}
\begin{figure}[!ht]
	\begin{subfigure}{0.33\textwidth}
		\includegraphics[width=\linewidth, height=4.5cm]{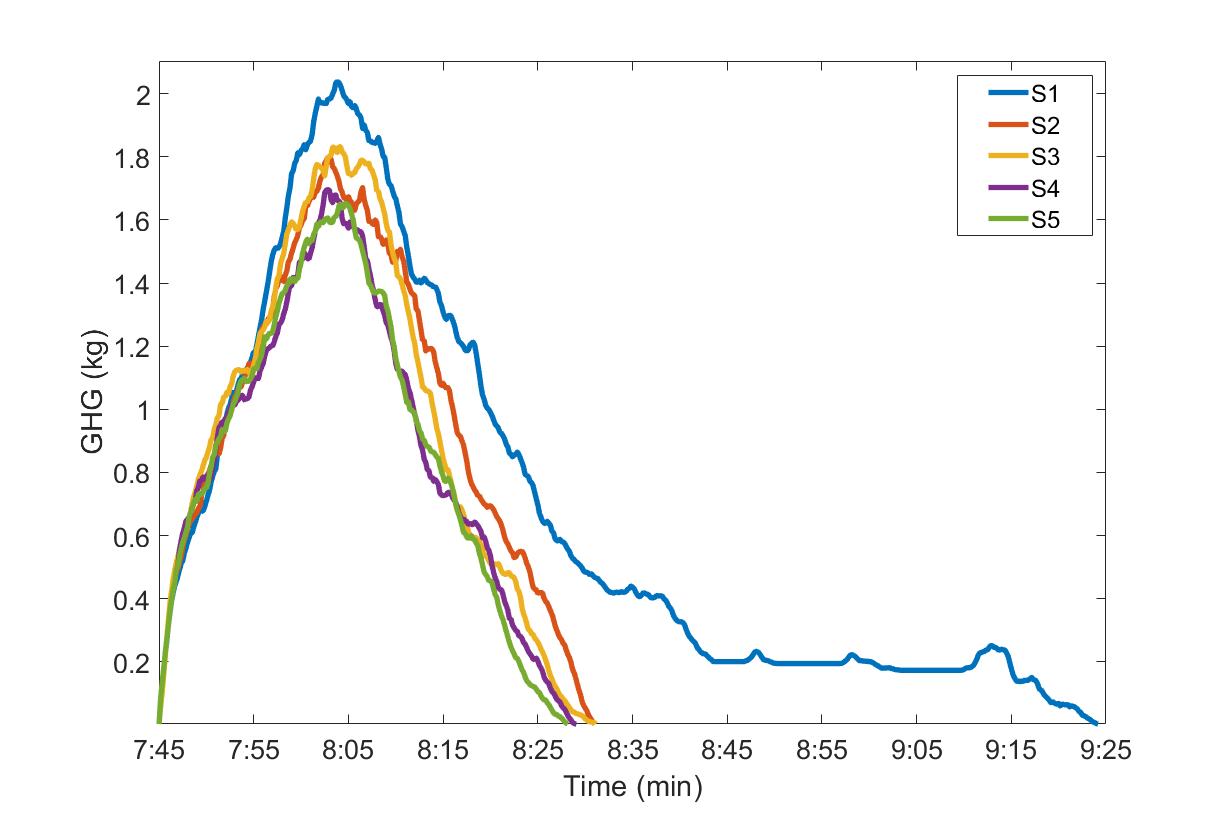} 
		\caption{GHG over time}
		\label{Network_GHG_over_time}
	\end{subfigure}%
	\begin{subfigure}{0.33\textwidth}
		\includegraphics[width=\linewidth, height=4.5cm]{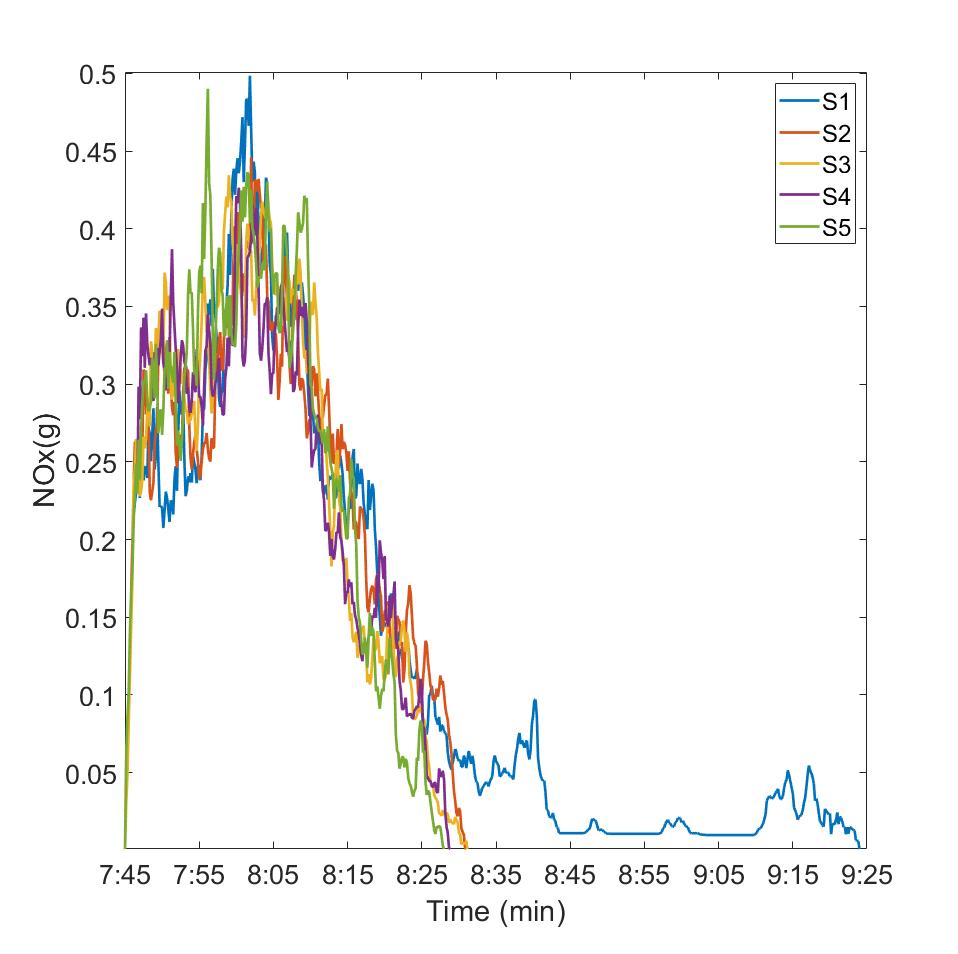}
		\caption{NOx over time}
		\label{Network_NOx_over_time}
	\end{subfigure}%
	\begin{subfigure}{0.33\textwidth}
		\includegraphics[width=\linewidth, height=4.5cm]{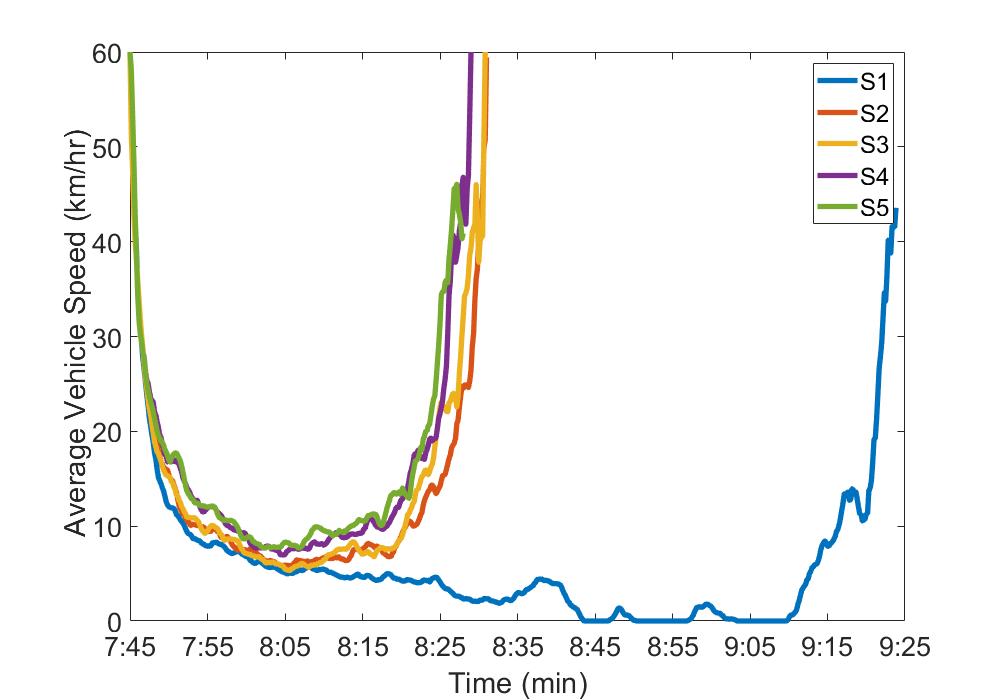}
		\caption{Average speed over time}
		\label{Network_average_speed_over_time}
	\end{subfigure}
	\caption{Network level GHG, NO$_X$, and average vehicle speed over time}
	\label{Network GHG, NOx_average speed over time}
\end{figure}

Figure \ref{Network GHG, NOx_average speed over time} illustrates results of network GHG, NO$_x$, and average speed over time for S1--S5, respectively. \textcolor{black}{It is worth mentioning that the GHG and NO$_x$ presented are the produced emission, the dispersion and exposure effects are not considered. As such when all vehicles leave the network at the end of the simulation the production levels for GHG and NO$_x$ are depicted as zero. In addition, both network and link average speed graphs presented in this section represent speed distribution while vehicles are still in the network}. There are four noticeable findings: 
\begin{itemize}
	\item[a.] The multi-objective objective function (Equation \ref{e:R2}) used in scenario S5 has the potential of reducing network throughput rate by unloading the network almost an hour faster than S1. This improvement is associated with the fact that the multi-objective function in S5 with the aid of network of communicating intelligent intersections using real-time traffic information tries to re-route vehicles not only to minimize travel time, but also to reduce idling time. Additionally, since the aim of Equation \ref{e:R2} is to also reduce the GHG emissions, it avoids unnecessary re-routing and formation of additional congestion by maintaining optimal speed level, while taking into account the VKT. In a way, this strategy is analogous to system optimal routing. It aims to reduce travel time, but takes into account the social cost of driving (e.g. queue at intersections and environmental pollution).
	\item[b.] In Figures \ref{Network_GHG_over_time}, \ref{Network_NOx_over_time}, \& \ref{Network_average_speed_over_time}, it can be observed that in terms of GHG, NO$_x$ and average network speed over time, during the first 10 to 15 minutes of the simulation (7:45-8:00), all routing strategies (S1--S5) perform almost the same. It is only around 8:00am when the demand loading is completed that the differences between routing scenarios become more apparent. As such around 8:00am is the most congested time on the network with most demand being loaded to the network. It can be noticed that the peak of congestion is roughly the same under all five scenarios. However, the severity of its effect on the network is not the same. Furthermore, it can be observed that after congestion formation, as we move from S1--S5 the ability of the network to recover from congestion increases with S5 outperforming the other scenarios. Up until the peak congestion period S1 is still capable of managing traffic efficiently. However, after congestion appears, S1's ability to react to changes on the network is reduced. In fact as can be observed in Figures \ref{Network_GHG_over_time} , \ref{Network_NOx_over_time}, \& \ref{Network_average_speed_over_time}, S1 takes almost an hour longer than the other four scenarios to recover from congestion, which in return results in longer travel times, higher GHG and NO$_x$ emission levels.
	\item[c.] From Figure \ref{Network GHG, NOx_average speed over time}, it can be observed that in the case of pre-trip dynamic routing under S1, there is a noticeable period when the average network travel speed is around 0km/hr, which is not the case in any other scenarios. This shows how without proper access to real-time traffic information vehicles may get stuck in queue in one part of the network, while not taking advantage of underutilized part of the network. Here, Front St. corridor has a speed limit of 80km/hr, as such in the case of S1 this is a frequently used corridor. However, a segment of Front St. that is between York and Bay St. (refer to Figure \ref{simplifiedtorontonetwork}), changes from 2 lanes street to one lane street and the speed drastically drops from 80km/hr to 40km/hr, creating the potential for gridlock at the intersection of Front St. with York St. and University Ave. This phenomenon also shows the importance of taking into consideration \textcolor{black}{idling penalty} at intersections when routing vehicles.
	\item[d.] By comparing Figures \ref{Network_GHG_over_time}, \ref{Network_NOx_over_time} and \ref{Network_average_speed_over_time}, it can be observed that at the beginning of the simulation for all five scenarios, the increase in NO$_x$ emission was more abrupt than increase in GHG and in fact, there were instances between 7:55am to 8:00am when the NO$_x$ was even higher under S5. \textcolor{black}{In the case of NO$_x$ the emission factor is more susceptible to higher speed, the faster vehicles travel more NO$_x$ they produce. Hence at the beginning of the simulation when there is less congestion and vehicles travel fast, NO$_x$ emissions level is higher. Lower speed produces less NO$_x$ and leads to vehicles spending more time on the network. Thus resulting in higher NO$_x$ emissions. As such, there needs to be a balance between how fast a vehicle should travel and how long it should stay in the network.} 
\end{itemize}

\begin{figure}[ht!]
	\centering
	\includegraphics[width=0.6\textwidth]{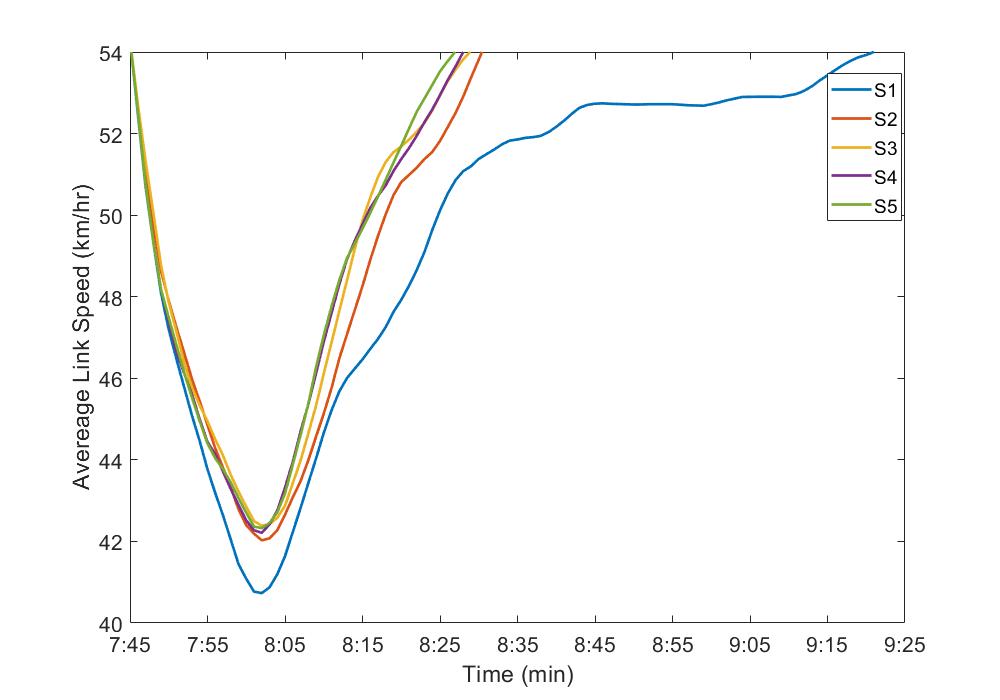}
	\caption{Average link speed over time.}\label{averagelinkspeed}
\end{figure}

Figure \ref{averagelinkspeed} presents the average link speed over time for the five scenarios. The difference between Figure \ref{averagelinkspeed} and \ref{Network_average_speed_over_time} is that Figure \ref{averagelinkspeed} is obtained using link level speeds where as Figure \ref{Network_average_speed_over_time} is obtained using individual vehicle speeds. Comparison between Figures \ref{averagelinkspeed} and \ref{Network_average_speed_over_time} reveals that at peak congestion period the average network link speed is noticeably higher than average vehicle speed on the network. It may be because even under real-time routing, it is possible that some links are avoided and they stay at their free flow speed or close to that while some links are saturated. When averaging over link speeds the links with free flow speeds are taken into account. Whereas, when averaging vehicle speeds those links are not counted. Hence, average vehicle speed is lower than actual average link speed. However, in Figure \ref{averagelinkspeed} it can be observed that as we move from S1 to S5 this difference between average link speed and average vehicle speed is minimized, because the network is being utilized better and vehicles are distributed over the entire network.

\subsection{Average GHG, \text{NO$_x$} and Speed per Trip}
\begin{figure}[!ht]
	\begin{subfigure}{0.33\textwidth}
		\includegraphics[width=\linewidth, height=4.5cm]{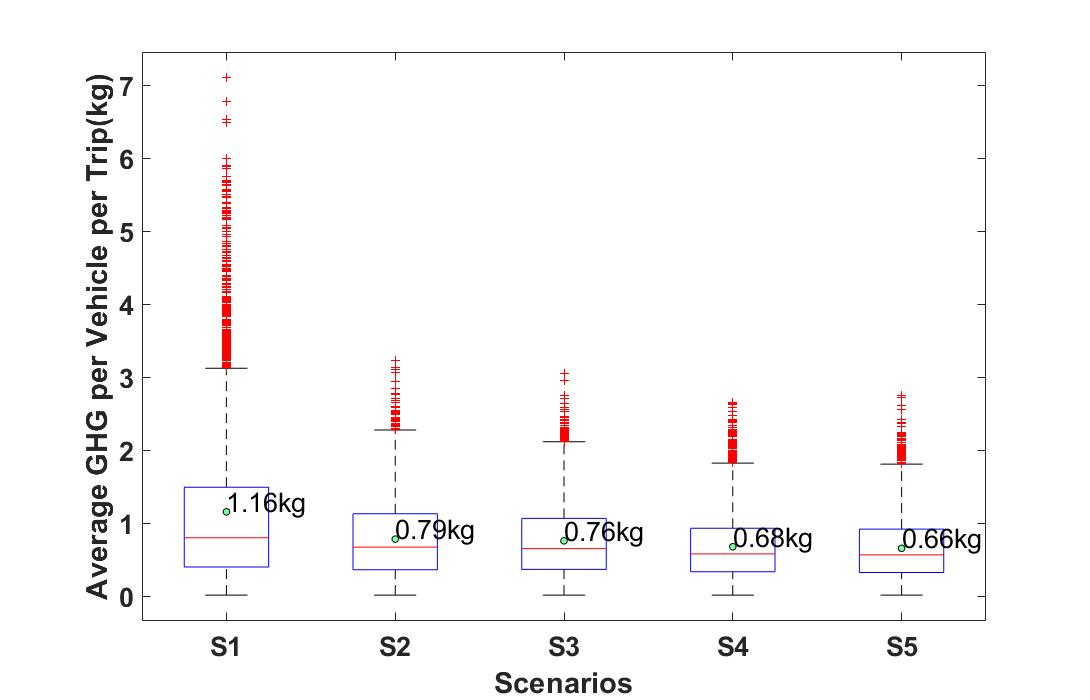} 
		\caption{Average GHG (kg)}
		\label{Network_GHG_boxplot}
	\end{subfigure}%
	\begin{subfigure}{0.33\textwidth}
		\includegraphics[width=\linewidth, height=4.5cm]{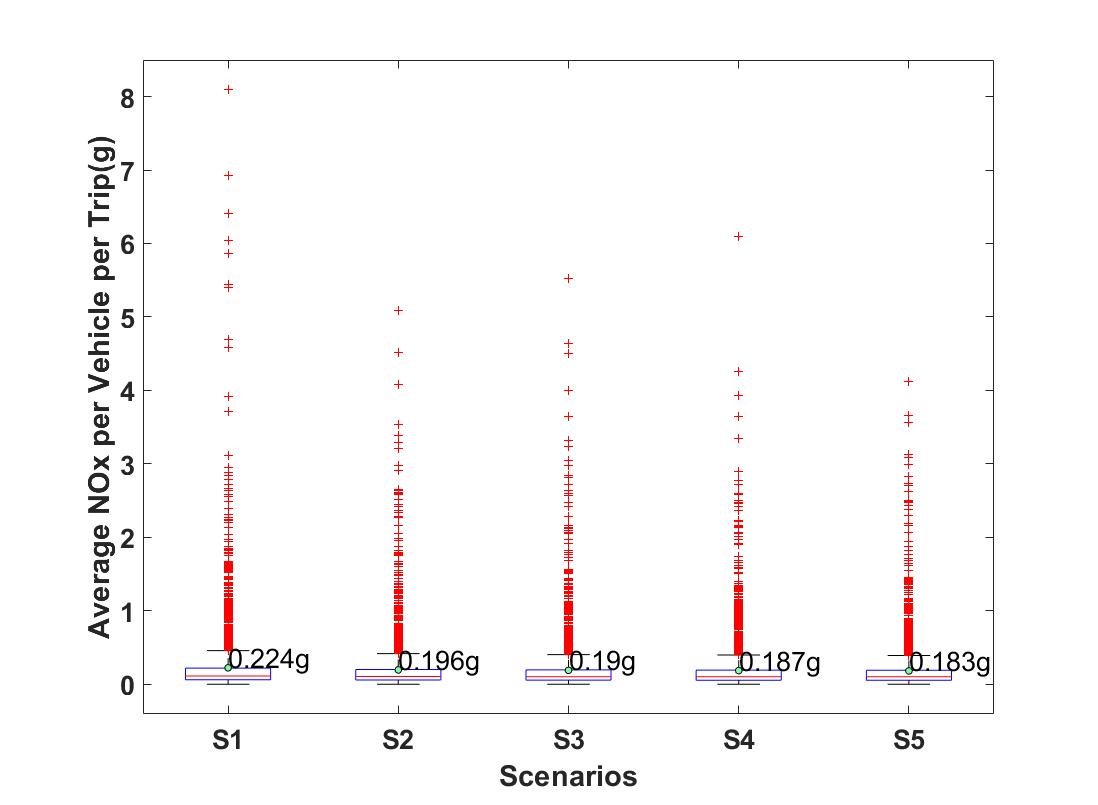}
		\caption{Average NOx (g)}
		\label{Network_NOx_boxplot}
	\end{subfigure}%
	\begin{subfigure}{0.33\textwidth}
		\includegraphics[width=\linewidth, height=4.5cm]{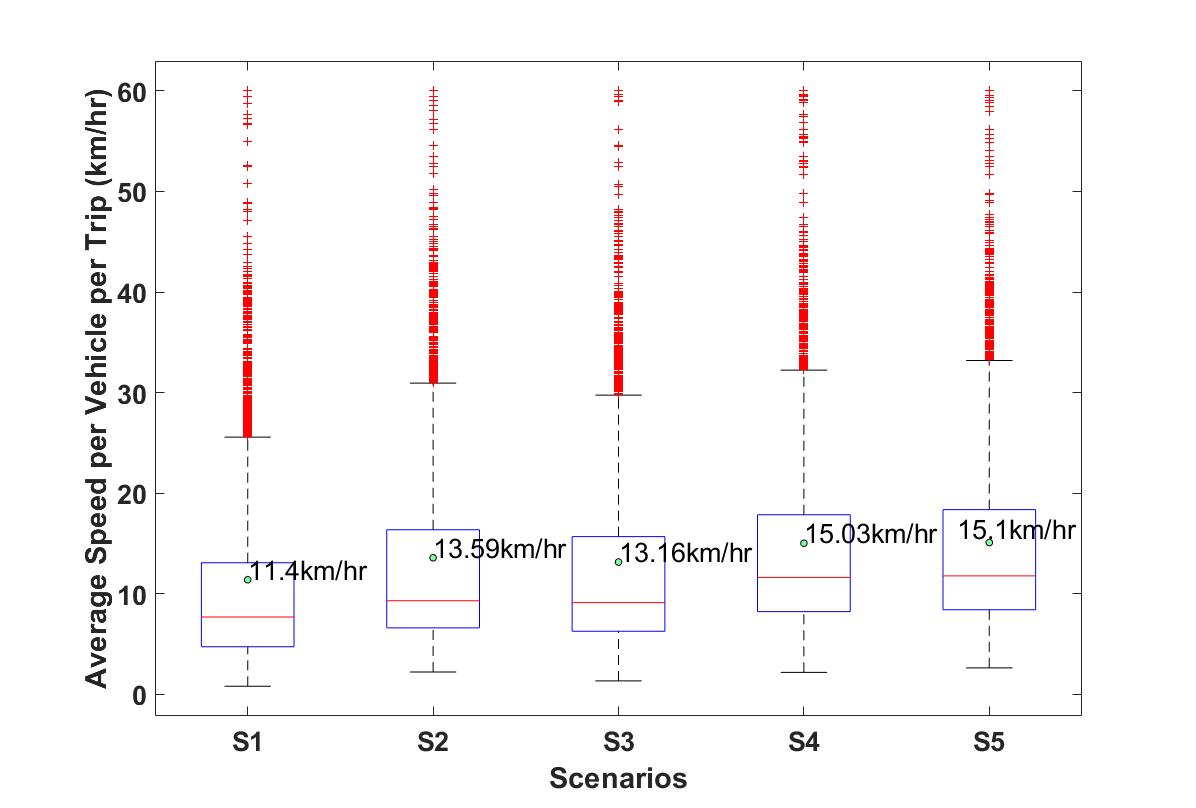}
		\caption{Average speed (km/hr)}
		\label{Network_average_spped_boxplot}
	\end{subfigure}
	\caption{Average GHG, NO$_X$, and vehicle speed}
	\label{average GHG, NOx and speed boxplot}
\end{figure}

The comparative analysis of average GHG, NO$_x$ and speed per trip is illustrated by Figure \ref{average GHG, NOx and speed boxplot}. Considering Figure \ref{Network_GHG_boxplot}, it can be seen that GHG emissions per trip have a higher variation under pre-trip dynamic routing (S1) in comparison to the other four scenarios (S2--S5), and as we go towards S5 this variation decreases. In scenario S1 with pre-trip dynamic routing, the trip-mean GHG emission is above the median suggesting a higher variation in the trips that emit GHG emissions higher than 3kg. Moving from S1 to S5, as the routing objective function takes more criteria into consideration, the variation in trip-level GHG emissions is reduced. In the case of S5 with routing strategy R2, the mean GHG emission is roughly the same as median, suggesting equal variation in both lower and upper 50\% of the data. In scenario S5, 75\% of the vehicles produced less than 2kg of GHG, which is a 50\% reduction in comparison to S1. This significant improvement is achieved by routing based on the real-time network-level traffic conditions to minimize travel time, delay and GHG emissions. Similar trend can be observed in average NO$_x$ emissions per trip shown in Figure \ref{Network_NOx_boxplot}, although not as noticeable as the case of GHG emissions. 

Figure \ref{Network_average_spped_boxplot} presents the boxplot comparative analysis of the five tested scenarios in terms of average trip speed. Under all the tested scenarios the average speed is higher than median, suggesting greater variation in upper 50\% of the data. However this variation is less in S1, which is due to the fact that vehicles are not re-routed en-route. Once they get stuck in congestion, they have to endure it till it dissipates. Whereas, in the S2--S5 vehicles are re-routed and they can experience different speed levels. Another interesting finding is that although the average speed per trip in scenario S3 under routing strategy R1 is approximately $15.5\%$ higher than S1. It is still lower than for the other routing scenarios S2, S4 and S5 where travel time is explicitly considered. In S3, intersections are constantly aware of the real-time emission levels, so they may be able to avoid congestion indirectly by avoiding the links with high GHG emissions. However, the speed in S3 compared to S2--S5 is lower due the quasi-convex relationship between GHG and speed. Overall, scenario 5 with routing objective R2 and E2ECAV as the routing scheme was able to outperform the other 4 scenarios and increase the average speed per trip by 32\% in comparison to scenario S1. 

\begin{figure}[!ht]
	\begin{subfigure}{0.5\textwidth}
		\includegraphics[width=\linewidth, height=4.5cm]{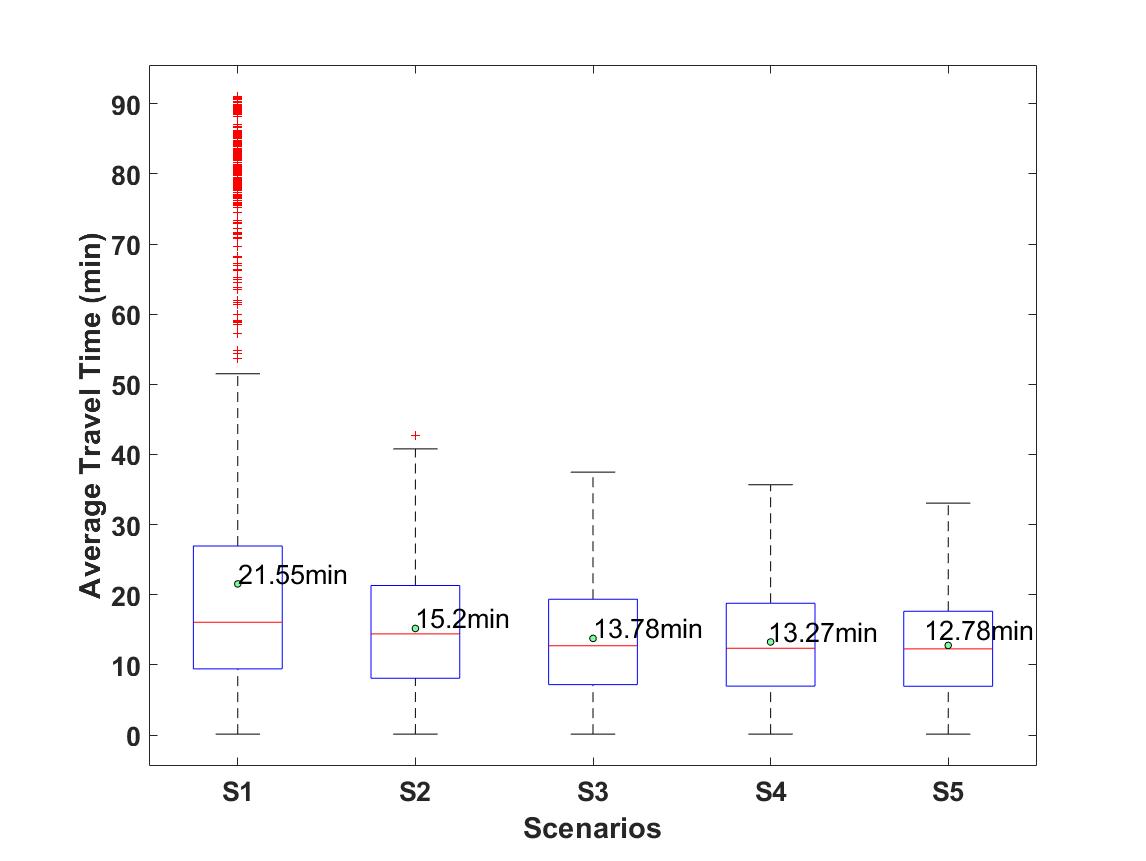} 
		\caption{Average travel time (min)}
		\label{average_TT_boxplot}
	\end{subfigure}%
	\hfill
	\begin{subfigure}{0.5\textwidth}
		\includegraphics[width=\linewidth, height=4.5cm]{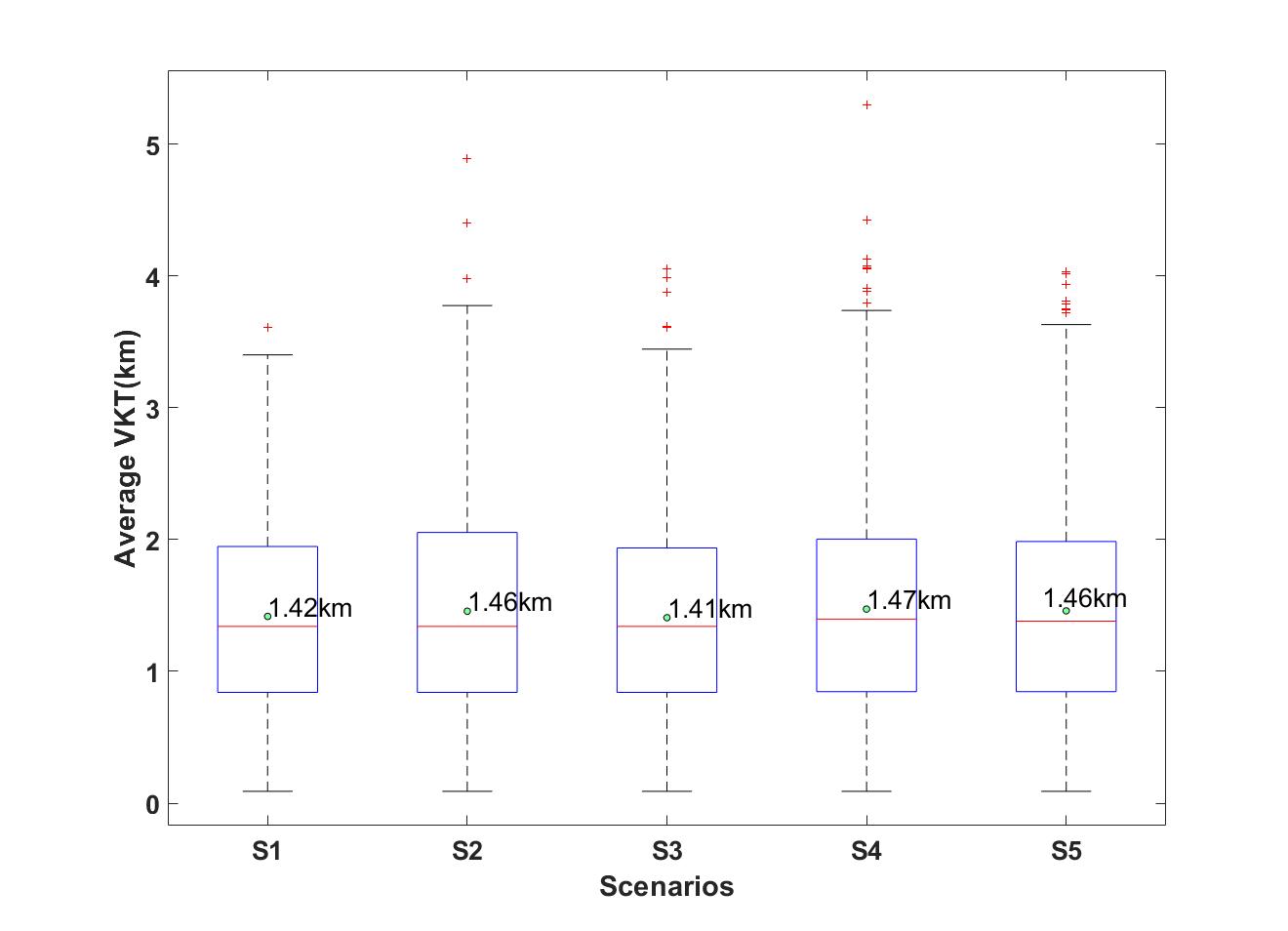}
		\caption{Average VKT}
		\label{Network_average_vkt_boxplot}
	\end{subfigure}
	\caption{Average travel time and average VKT}
	\label{average TT and VKT boxplot}
\end{figure}

Figure \ref{average TT and VKT boxplot} presents the boxplot analysis for average travel time (min) and average VKT on the network. It can be seen that under all 5 scenarios both the averages are above the medians, suggesting greater degree of variation in the upper 50\% of the data. However as we move from S1 to S5, this variation decreases. In Figure \ref{average_TT_boxplot}, scenario 5 has the least amount of variation, with mean travel time being close to median, this is due to the fact that E2ECAV routing scheme used here is better at avoiding both current congestion as well as avoiding new ones. This further demonstrates that routing strategy R2 used in scenario S5 is able to manage congestion better by considering travel time, delay and GHG emission while routing vehicles. In Figure \ref{Network_average_vkt_boxplot}, it is noticeable that scenario S3 had the least variation with mean VKT being close to median. The reason for this is that in S3 the system is more sensitive to GHG emissions than travel time---as such vehicles are only re-routed when there is a significant change in GHG emissions. In case of $TT$ based routing (S2, S4, and S5), when there is congestion, vehicles are re-routed immediately to a faster moving links. Thus, leading to higher re-routing rate and variation in VKT. By comparing Figure \ref{Network_average_spped_boxplot}, \ref{average_TT_boxplot} and \ref{Network_average_vkt_boxplot}, it can be seen that although scenario S3 has a lower average speed than S2, S4, and S5, it still has a lower average travel time due to the shorter trips resulting from routing strategy R1.

\textcolor{black}{Statistical significance $t-$test was performed to test the significance of proposed multi-objective algorithms under different scenarios in terms of $TT$, $VKT$, speed, $GHG$ and $NO_x$. Table \ref{sigtest} presents the results for the statistical significance test. Table \ref{sigtest} shows that each of the tested routing scenarios S2-S5 resulted in significant changes in terms of TT, speed and GHG not only in comparison to S1, but also in comparison to each other. These results are inline with the findings presented in previous sections and further emphasize the necessity for real time routing, along with inclusion of environmental pollution and delay penalty. In terms of VKT it can be seen that with the exception of S3 all the other scenarios resulted in significant change in comparison to S1. This noticeable change was due to ability of S2, S4 and S5 en-route routing and congestion avoidance, whereas in case of S3 en-route routing was more restrictive in order to avoid substantial increase in VKT. In general, wherever the travel time was the dominant variable in the objective function a significant change in VKT is apparent. From Table \ref{sigtest} it can be observed that S2-S5 were able to significantly affect NO$_x$ in comparison to S1, this is in line with earlier results and is due to the fact that S2-S5, as shown in Figure \ref{TT}, were able to reduce travel time substantially compared to S1. The less time vehicles spend in less NO$_x$ they will produce. In the scenarios where en-route routing was available (S2-S5), it can be seen that when travel time, idling penalty and GHG were considered explicitly and simultaneously it resulted in most significant improvement in NO$_x$. }
\begin{table}[!ht]
	\caption{\textcolor{black}{Statistical significance \textit{t}-test}}
	\label{sigtest}
	\centering
	\small
	\begin{tabular}{|c|c|c|c|c|c|c|c|c|c|c| }
		\cline{1-11}
		\multicolumn{1}{|l|}{} & \multicolumn{2}{|c|}{\textbf{TT}}                                      & \multicolumn{2}{|c|}{\textbf{VKT}}                                     & \multicolumn{2}{|c|}{\textbf{Speed}}                                   & \multicolumn{2}{|c|}{\textbf{GHG}}                                     & \multicolumn{2}{|c|}{\textbf{NO$_x$}}                                     \\
		\cline{1-11}
		\textbf{Scenarios}   & \multicolumn{1}{c}{\textbf{t-stat}} & \multicolumn{1}{|c|}{\textbf{p}} & \multicolumn{1}{c}{\textbf{t-stat}} & \multicolumn{1}{|c|}{\textbf{p}} & \multicolumn{1}{c}{\textbf{t-stat}} & \multicolumn{1}{|c|}{\textbf{p}} & \multicolumn{1}{c}{\textbf{t-stat}} & \multicolumn{1}{|c|}{\textbf{p}} & \multicolumn{1}{c}{\textbf{t-stat}} & \multicolumn{1}{|c|}{\textbf{p}} \\
		\cline{1-11}
		\textbf{S1 vs. S2}   & 17.653                              & 0.000                          & -2.152                              & 0.031                          & -8.365                              & 0.000                          & 17.877                              & 0.000                          & 3.050                               & 0.002                          \\
		\hline
		\textbf{S1 vs. S3}   & 21.927                              & 0.000                          & \textcolor{blue}{0.568}                               & \textcolor{blue}{0.570}                          & -6.829                              & 0.000                          & 19.292                              & 0.000                          & 3.692                               & 0.000                          \\
		\hline
		\textbf{S1 vs. S4}   & 23.527                              & 0.000                          & -3.014                              & 0.003                          & -14.238                             & 0.000                          & 23.578                              & 0.000                          & 4.124                               & 0.000                          \\
		\hline
		\textbf{S1 vs. S5}   & 25.164                              & 0.000                          & -2.299                              & 0.022                          & -14.835                             & 0.000                          & 24.774                              & 0.000                          & 4.706                               & 0.000                          \\
		\hline
		\textbf{S2 vs. S3}   & 7.047                               & 0.000                          & 2.708                               & 0.007                          & 1.665                               & 0.096                          & 1.916                               & 0.055                          & \textcolor{blue}{0.754}                               & \textcolor{blue}{0.451}                          \\
		\hline
		\textbf{S2 vs. S4}   & 9.791                               & 0.000                          & \textcolor{blue}{-0.863}                              & \textcolor{blue}{0.388}                          & -5.555                              & 0.000                          & 8.713                               & 0.000                          & \textcolor{blue}{1.220}                               & \textcolor{blue}{0.222}                          \\
		\hline
		\textbf{S2 vs. S5}   & 12.670                              & 0.000                          & \textcolor{blue}{-0.149}                              & \textcolor{blue}{0.882}                          & -5.957                              & 0.000                          & 10.642                              & 0.000                          & 1.845                               & 0.065                          \\
		\hline
		\textbf{S3 vs. S4}   & 2.738                               & 0.006                          & -3.564                              & 0.000                          & -7.370                              & 0.000                          & 7.004                               & 0.000                          & \textcolor{blue}{0.451}                               & \textcolor{blue}{0.652}                          \\
		\hline
		\textbf{S3 vs. S5}   & 5.539                               & 0.000                          & -2.853                              & 0.004                          & -7.822                              & 0.000                          & 8.989                               & 0.000                          & 1.052                               & 0.293                          \\
		\hline
		\textbf{S4 vs. S5}   & 2.769                               & 0.006                          & \textcolor{blue}{0.714}                               & \textcolor{blue}{0.475}                          & \textcolor{blue}{-0.296}                              & \textcolor{blue}{0.767}                          & 1.926                               & 0.054                          & \textcolor{blue}{0.599}                               & \textcolor{blue}{0.549}         \\       
		\hline
	\end{tabular}
\end{table}

\subsection{Heatmaps: Speed, current vs jam densities, and GHG emissions}
Figures \ref{speed, density and GHG heatmap of hdv}, \ref{speed, density and GHG heatmap of TT}, \ref{speed, density and GHG heatmap of GHG}, \ref{speed, density and GHG heatmap of TTstar} and \ref{speed, density and GHG heatmap of TTstarplusGHG} present heatmaps of speed, current vs jam densities $\displaystyle \Bigg(\frac{\rho_{(\Delta_j,l)}}{\rho_{jam}}\Bigg)$, and GHG for scenarios S1--S5 at 8:20am. The following noticeable differences can be observed that are aligned with the previous results:\par
\begin{figure}[!ht]
	\begin{subfigure}{0.33\textwidth}
		\includegraphics[width=\linewidth=]{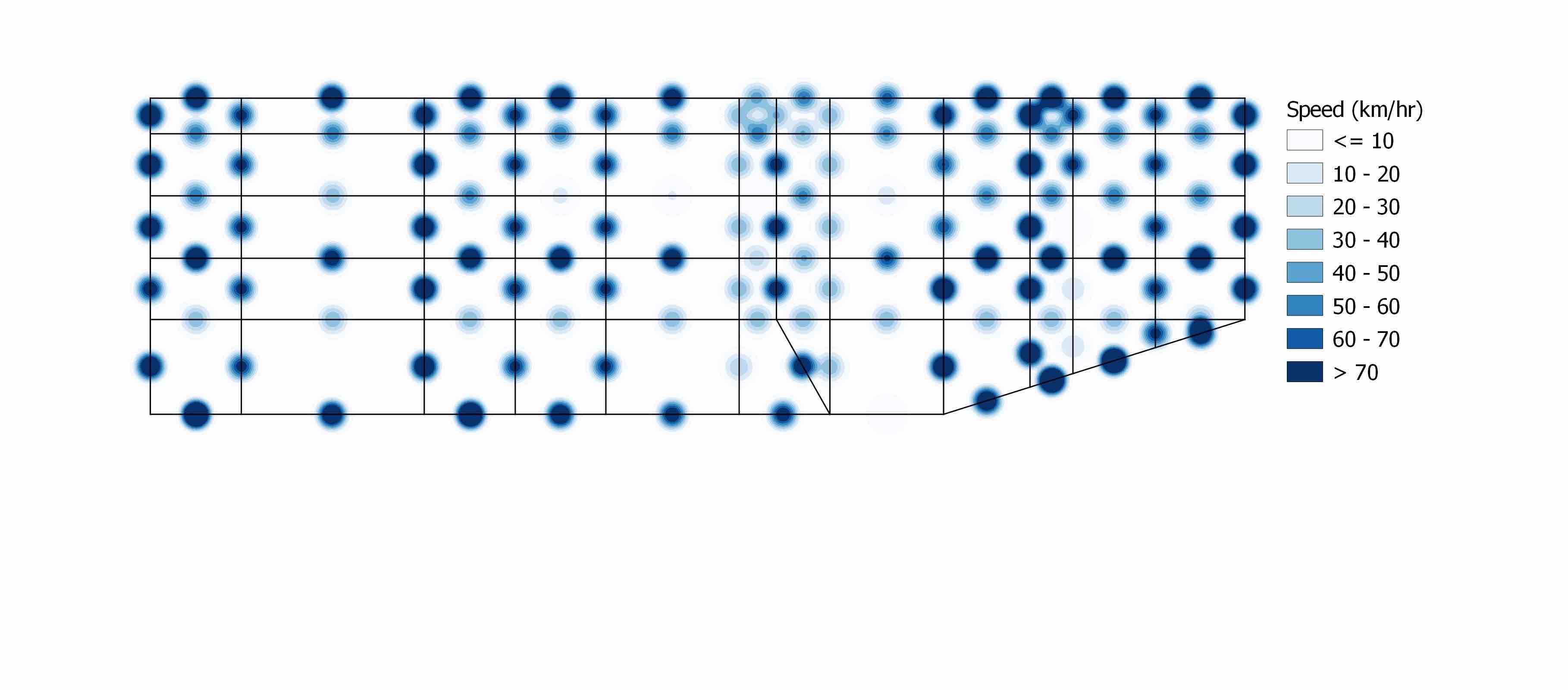} 
		\caption{Speed heatmap for S1}
		\label{speedheatmaphdv}
	\end{subfigure}%
	\begin{subfigure}{0.33\textwidth}
		\includegraphics[width=\linewidth=]{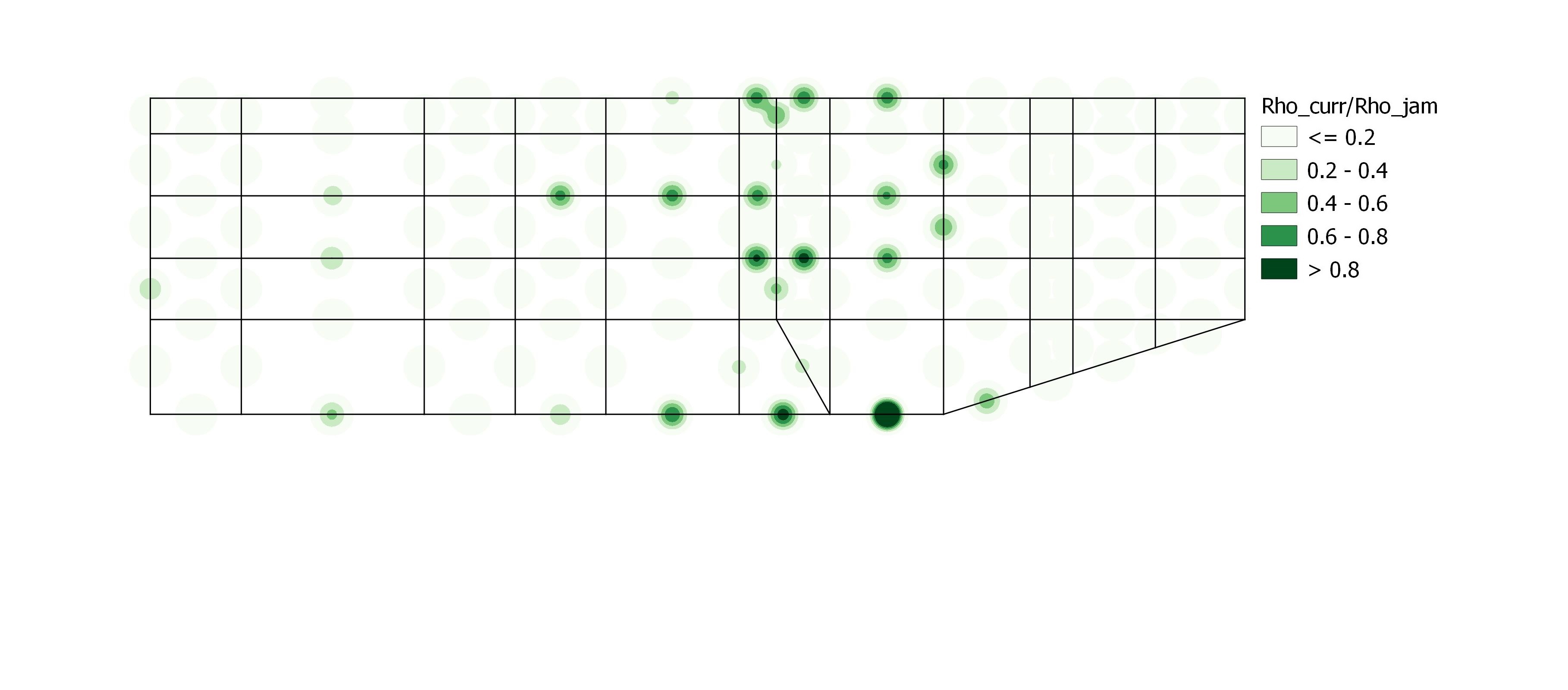}
		\caption{Density heatmap of S1}
		\label{Densityheatmapofhdv}
	\end{subfigure}%
	\begin{subfigure}{0.33\textwidth}
		\includegraphics[width=\linewidth]{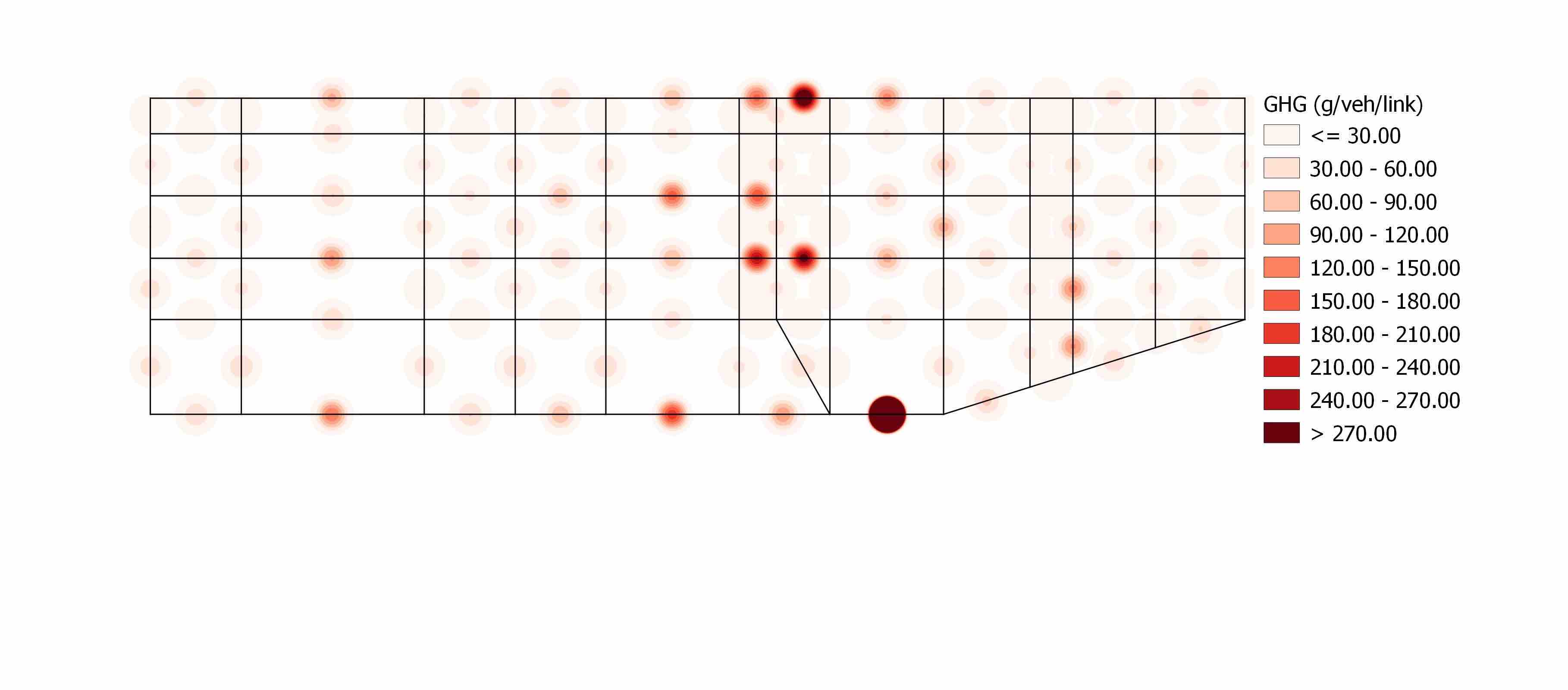}
		\caption{GHG heatmap of S1}
		\label{GHGheatmapofhdv}
	\end{subfigure}
	\caption{Speed, density and GHG heatmaps of S1}
	\label{speed, density and GHG heatmap of hdv}
\end{figure}

 
\begin{itemize}
	\item[a.] By comparing the heatmap of GHG emissions from S1 (Figure \ref{GHGheatmapofhdv}) with the other scenarios, it can be observed that S1 produces higher GHG emissions, which is concentrated at or around the intersections Front at University/York, King St at University, Adelaide at University, and Queen at University. This is due to high $\frac{\rho_{(\Delta_j,l)}}{\rho_{jam}}$ (\ref{Densityheatmapofhdv}) and subsequently lower speed on the associated links (\ref{speedheatmaphdv}). Those links are frequently used because they are the major arterials of the network. For example, Front St is a bidirectional street with two lanes and speed of 80km/hr, University Avenue is a bi-directional street with 3 lanes and speed of 60km/hr, and Adelaide is one way street (W-E) with 3 lanes. Due to their attractiveness, the above mentioned corridors are highly utilized in S1---as a result, they get congested faster. Once congested the system does not have the ability to recover quickly due to the lack of en-route routing. It can be seen that in scenarios S2--S5, with the implementation of E2ECAV based routing, there is a significant reduction in the GHG emissions and congestion level, along with an increase in speed.
	\item[b.] By comparing S3 (Figure \ref{speed, density and GHG heatmap of GHG}) with the results from other the four scenarios (S1, S2, S4, and S5), it can be observed that under routing strategy R1, Front St is not utilized as much as the other scenarios, the reason behind this is that speed limit of Front street is 80km/hr and this speed produces high level of emissions, as such it is not frequently used. This further verifies the quasi-convex relationship between speed and GHG emissions.
	\item[c.] In Figure \ref{speed, density and GHG heatmap of hdv}, \ref{speed, density and GHG heatmap of TT}, \ref{speed, density and GHG heatmap of GHG}, \ref{speed, density and GHG heatmap of TTstar} and \ref{speed, density and GHG heatmap of TTstarplusGHG}, it can be observed that the lower part of Victoria St. between King and Front St. under all five tested scenarios has high GHG emissions, despite of the fact that it has a very low $\frac{\rho_{(\Delta_j,l)}}{\rho_{jam}}$ ($<$ 0.2). The reason behind this phenomenon is that the section of Victoria St. changes from two lanes to one lane and speed drops from 40km/hr to 10km/hr, which is quite low and produces more emissions.
	\item[d.] Front street is highly used in S1, S4, and S5. However, unlike S1, links in S4 and S5 are less congested, with lower GHG emissions, and higher speeds. This improvement in congestion level and GHG emissions is due to the fact that under S4 and S5, intersections have access to delay information, real-time travel time and GHG emissions. They route vehicles to Front St. in moderation, avoiding creation of congestion and gridlock. Even if the congestion occurs, they are able to dissipate it faster by re-routing vehicles en-route.
	\item[e.] In S2 (Figure \ref{speed, density and GHG heatmap of TT}), it can be observed that $TT$ only based routing (Equation \ref{e:TT}) with its en-route routing capability and access to real-time traffic information was able to reduce congestion and GHG emissions in comparison to S1. However, the improvement is still not significant under S5. This further emphasizes the fact that minimizing only the travel time is not enough to decrease GHG emissions. A good routing strategy requires inclusion of travel time, \textcolor{black}{idling penalty}, and emissions. Overall, in comparison to routing strategies used in S1--S4, the routing strategy in S5 is able to distribute traffic more evenly throughout the network and it ensures higher speed, lower congestion and lower GHG emissions level. 
\end{itemize}

\newpage
\begin{figure}[!ht]
	\begin{subfigure}{0.33\textwidth}
		\includegraphics[width=\linewidth]{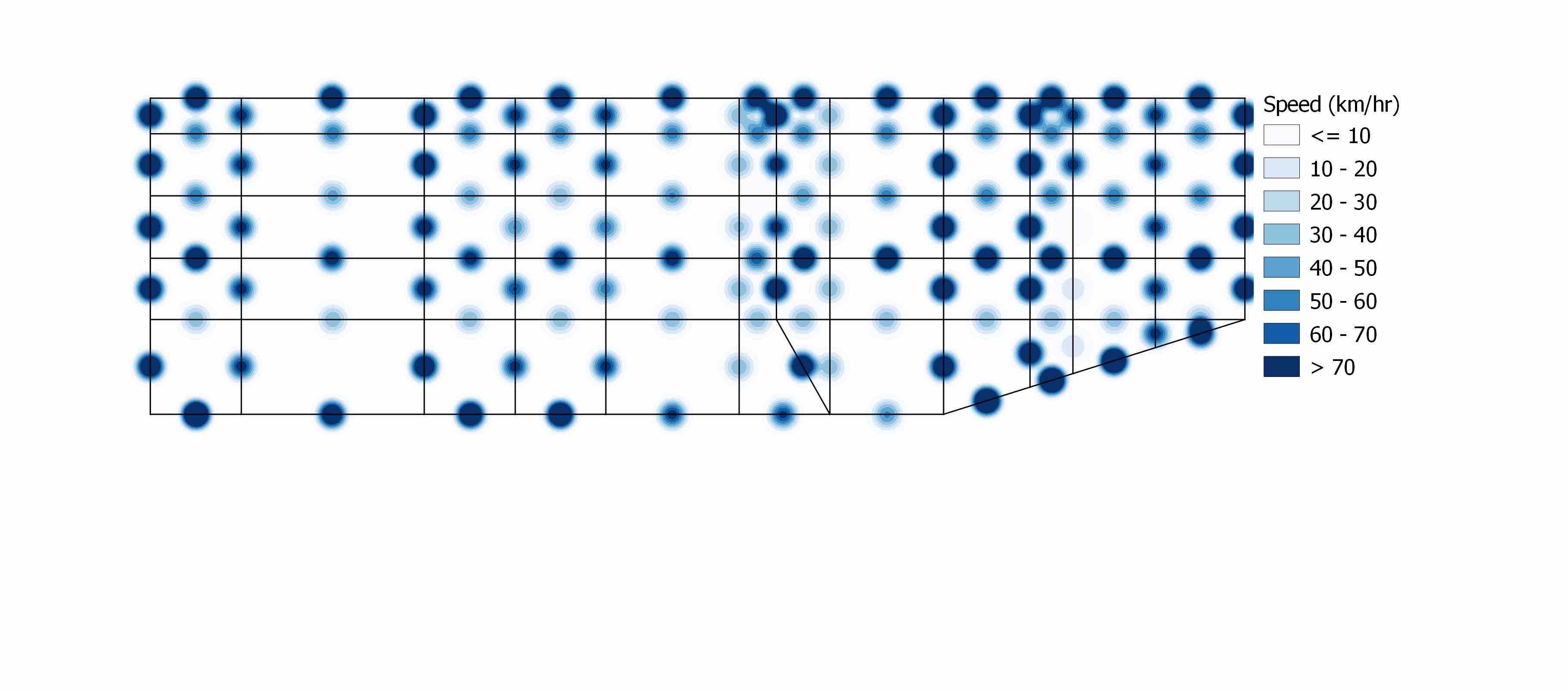} 
		\caption{Speed heatmap for S2}
		\label{speedheatmapTT}
	\end{subfigure}%
	\begin{subfigure}{0.33\textwidth}
		\includegraphics[width=\linewidth]{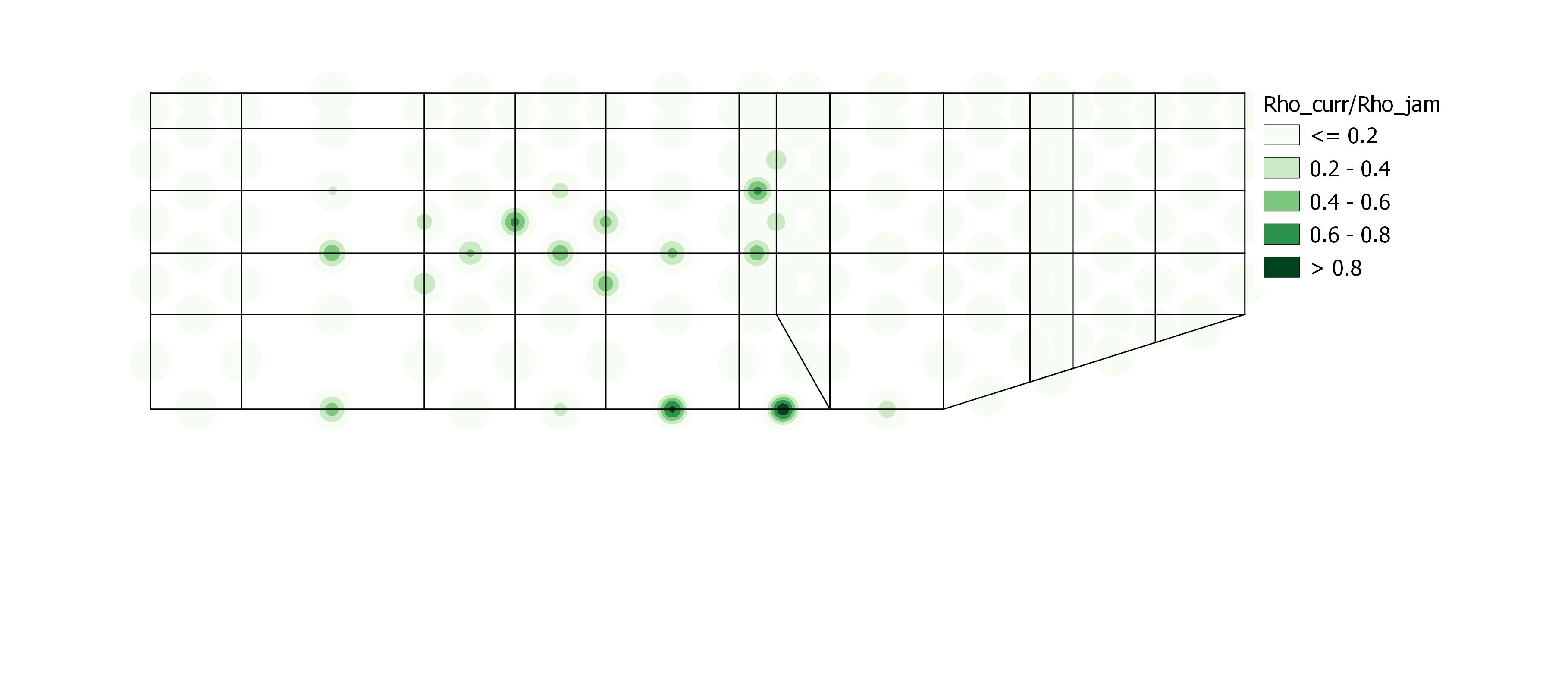}
		\caption{Density heatmap of S2}
		\label{DensityheatmapofTT}
	\end{subfigure}%
	\begin{subfigure}{0.33\textwidth}
		\includegraphics[width=\linewidth]{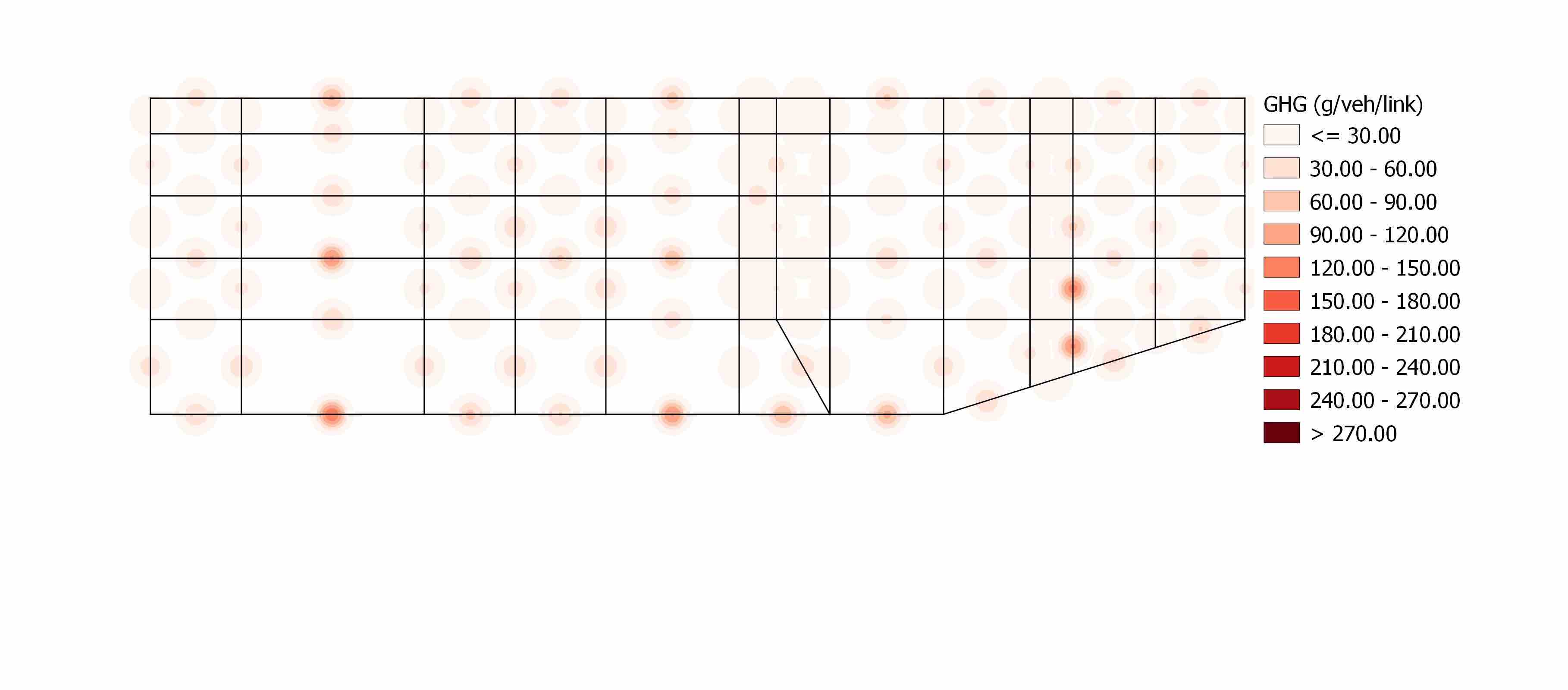}
		\caption{GHG heatmap of S2}
		\label{GHGheatmapofTT}
	\end{subfigure}
	\caption{Speed, density and GHG heatmaps of S2}
	\label{speed, density and GHG heatmap of TT}
\end{figure}
\begin{figure}[!ht]
	\begin{subfigure}{0.33\textwidth}
		\includegraphics[width=\linewidth]{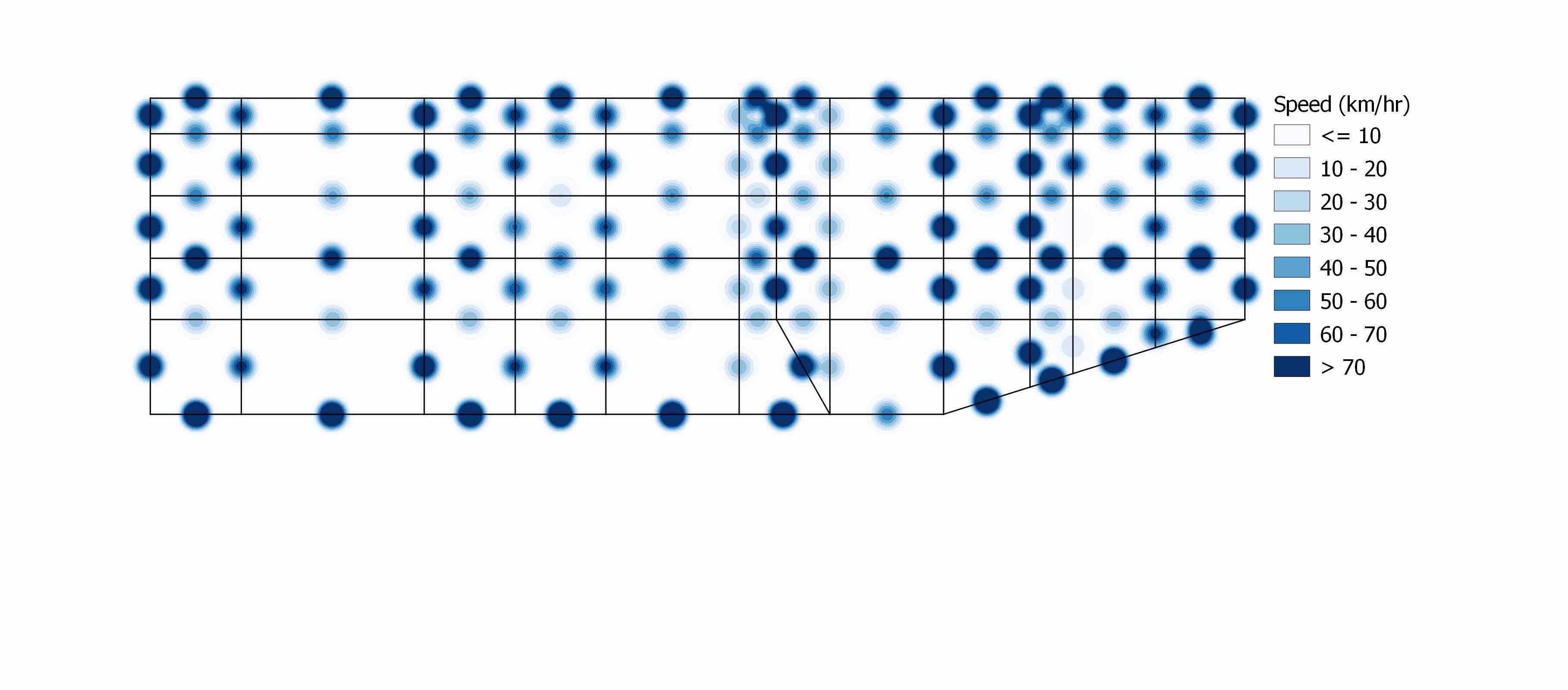} 
		\caption{Speed heatmap for S3}
		\label{speedheatmapGHG}
	\end{subfigure}%
	\begin{subfigure}{0.33\textwidth}
		\includegraphics[width=\linewidth]{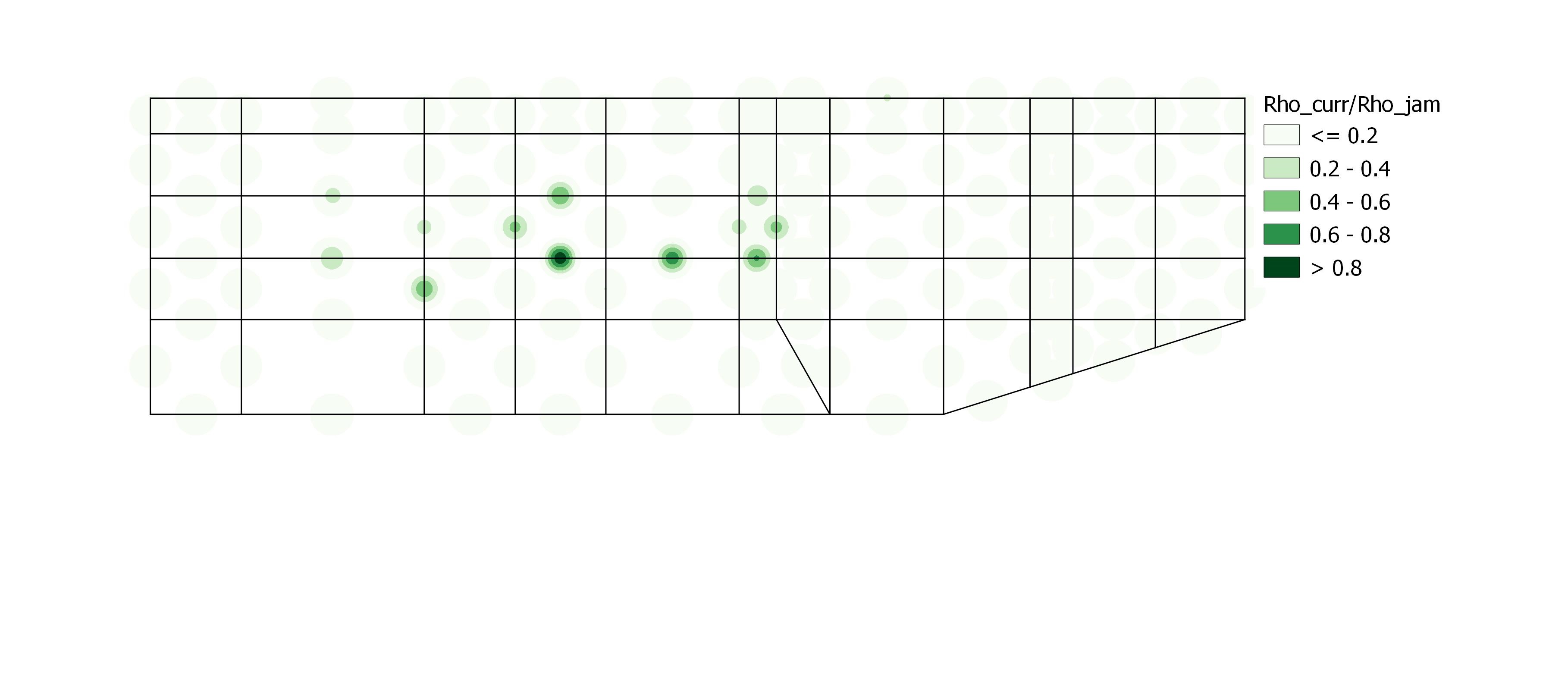}
		\caption{Density heatmap of S3}
		\label{DensityheatmapofGHG}
	\end{subfigure}%
	\begin{subfigure}{0.33\textwidth}
		\includegraphics[width=\linewidth]{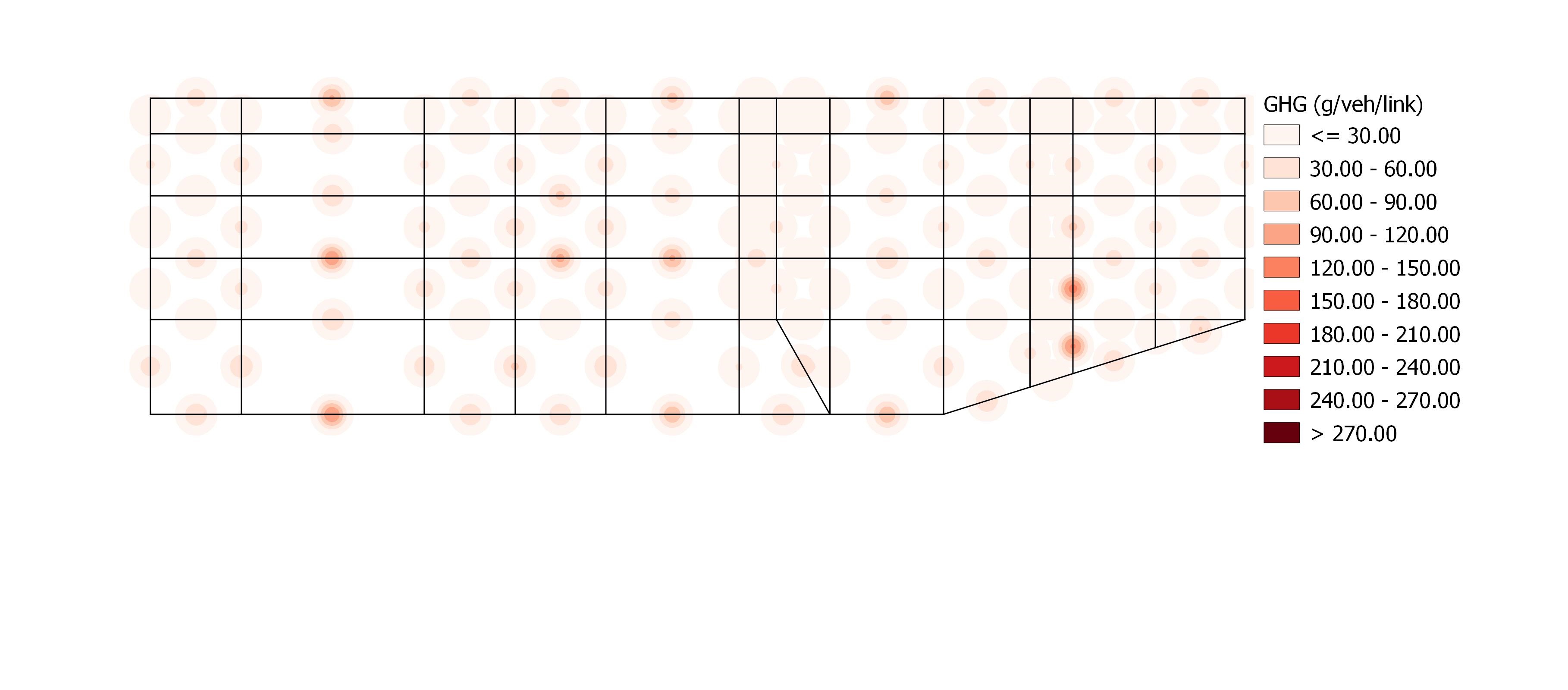}
		\caption{GHG heatmap of S3}
		\label{GHGheatmapofGHG}
	\end{subfigure}
	\caption{Speed, density and GHG heatmaps of S3}
	\label{speed, density and GHG heatmap of GHG}
\end{figure}
\begin{figure}[!ht]
	\begin{subfigure}{0.33\textwidth}
		\includegraphics[width=\linewidth]{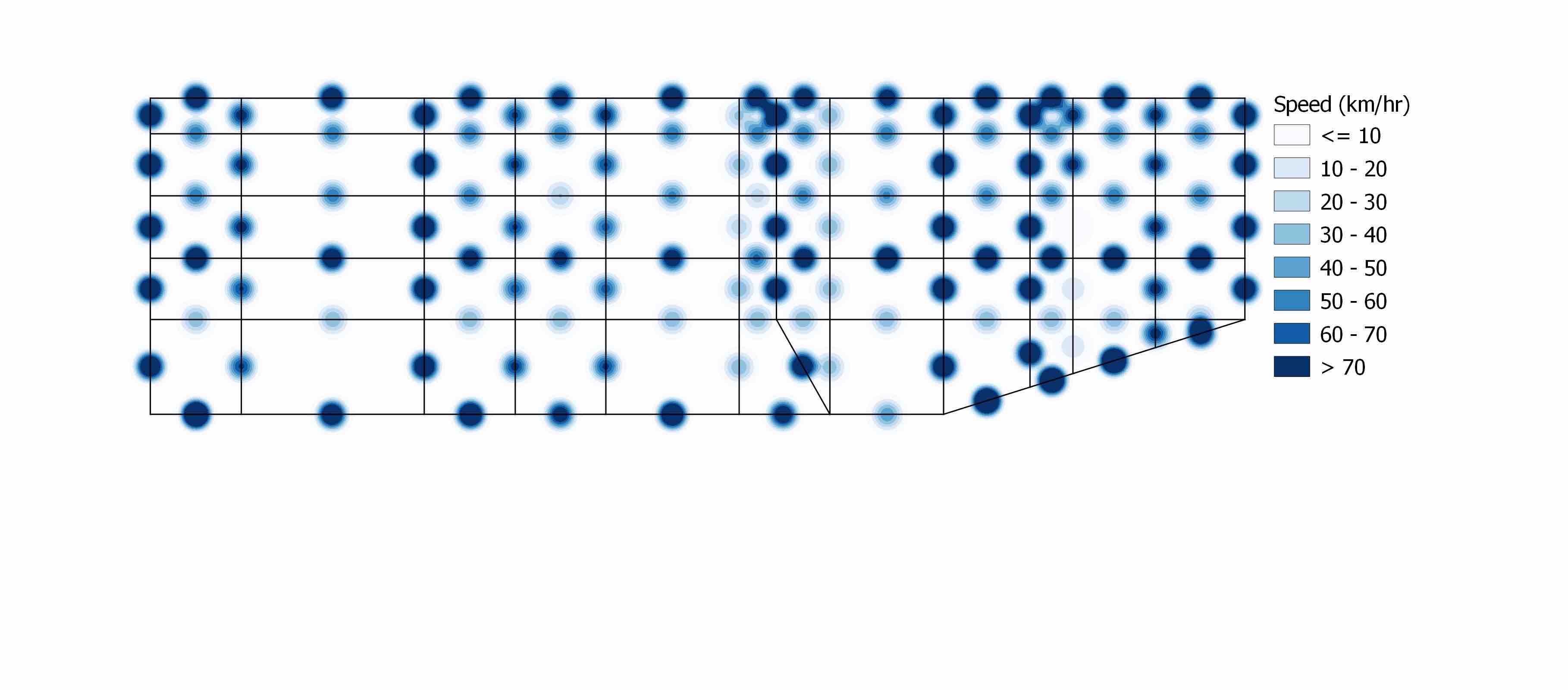} 
		\caption{Speed heatmap for S4}
		\label{speedheatmapTTstar}
	\end{subfigure}%
	\begin{subfigure}{0.33\textwidth}
		\includegraphics[width=\linewidth]{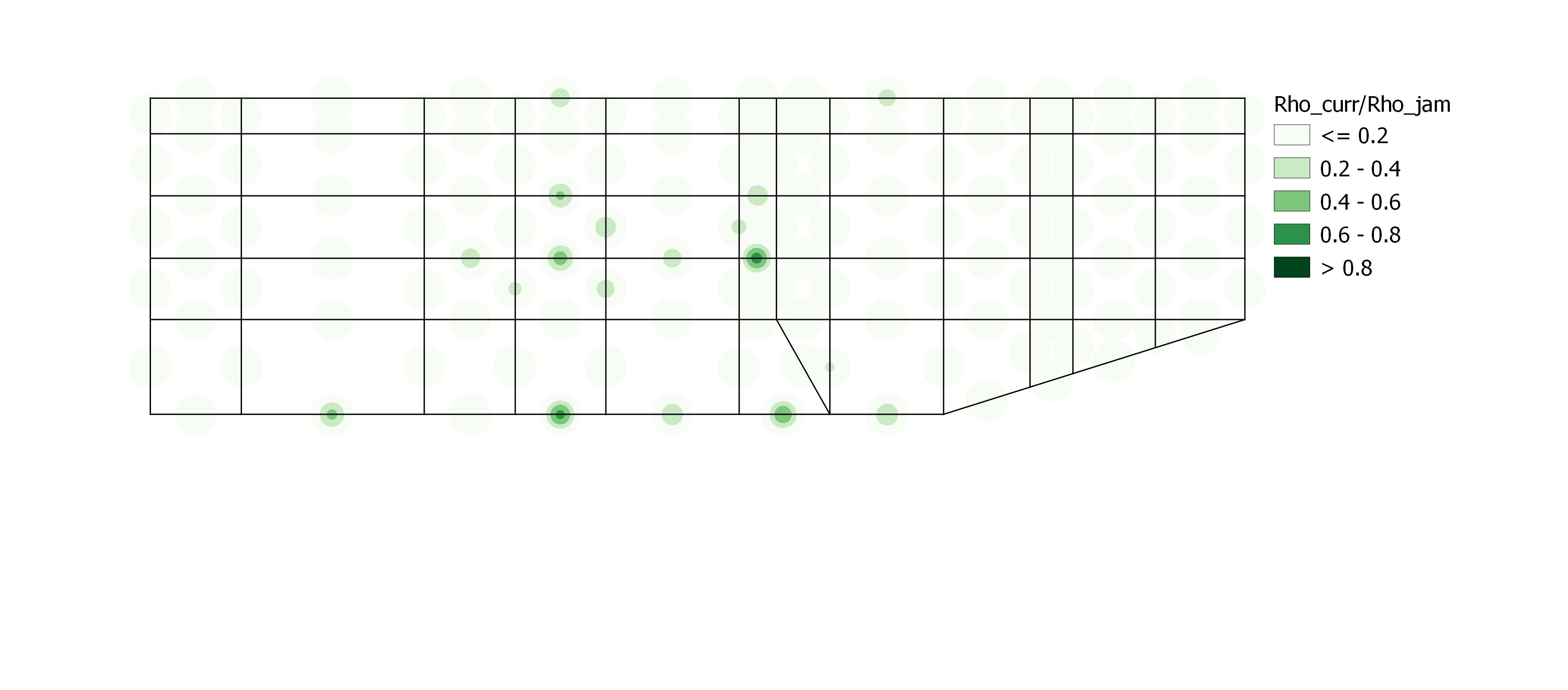}
		\caption{Density heatmap of S4}
		\label{DensityheatmapofTTstar}
	\end{subfigure}%
	\begin{subfigure}{0.33\textwidth}
		\includegraphics[width=\linewidth]{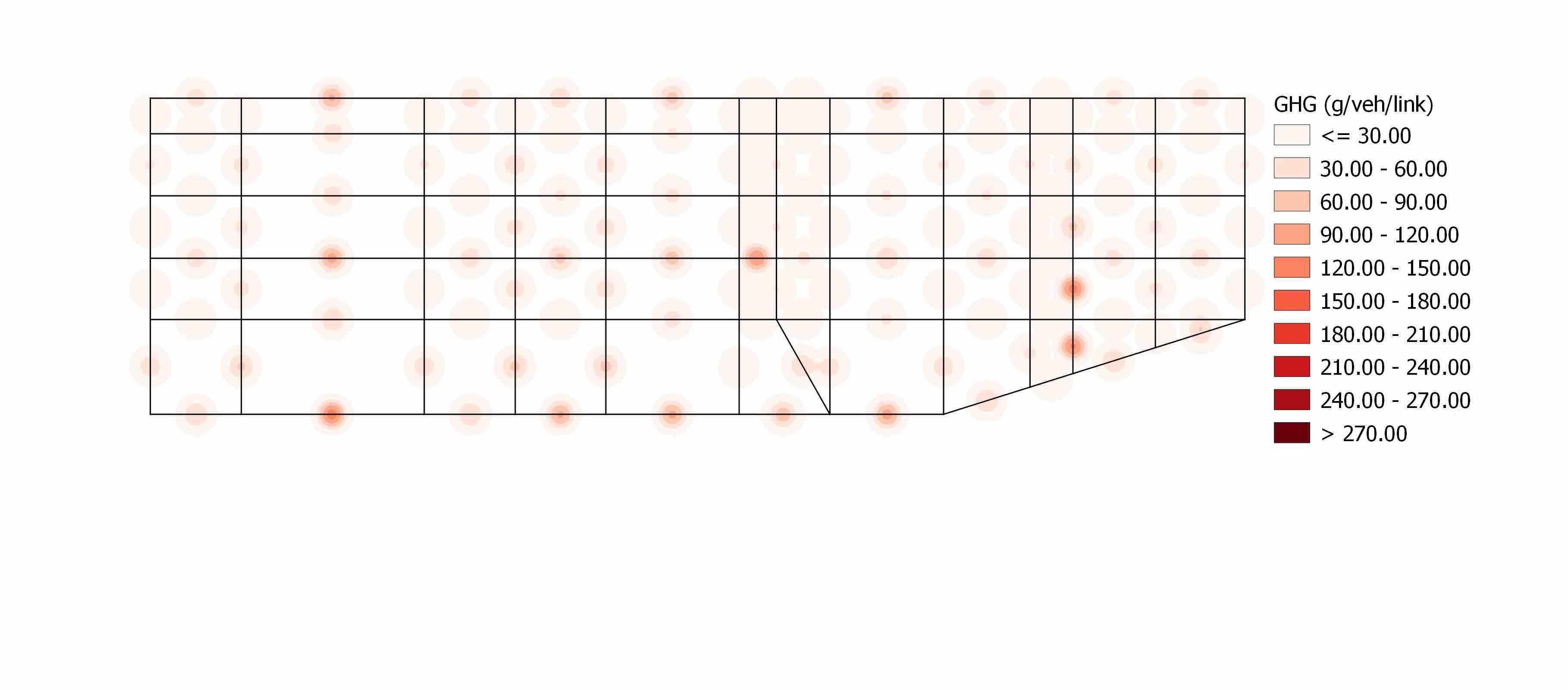}
		\caption{GHG heatmap of S4}
		\label{TTstarheatmapofGHG}
	\end{subfigure}
	\caption{Speed, density and GHG heatmaps of S4}
	\label{speed, density and GHG heatmap of TTstar}
\end{figure}
\begin{figure}[!ht]
	\begin{subfigure}{0.33\textwidth}
		\includegraphics[width=\textwidth]{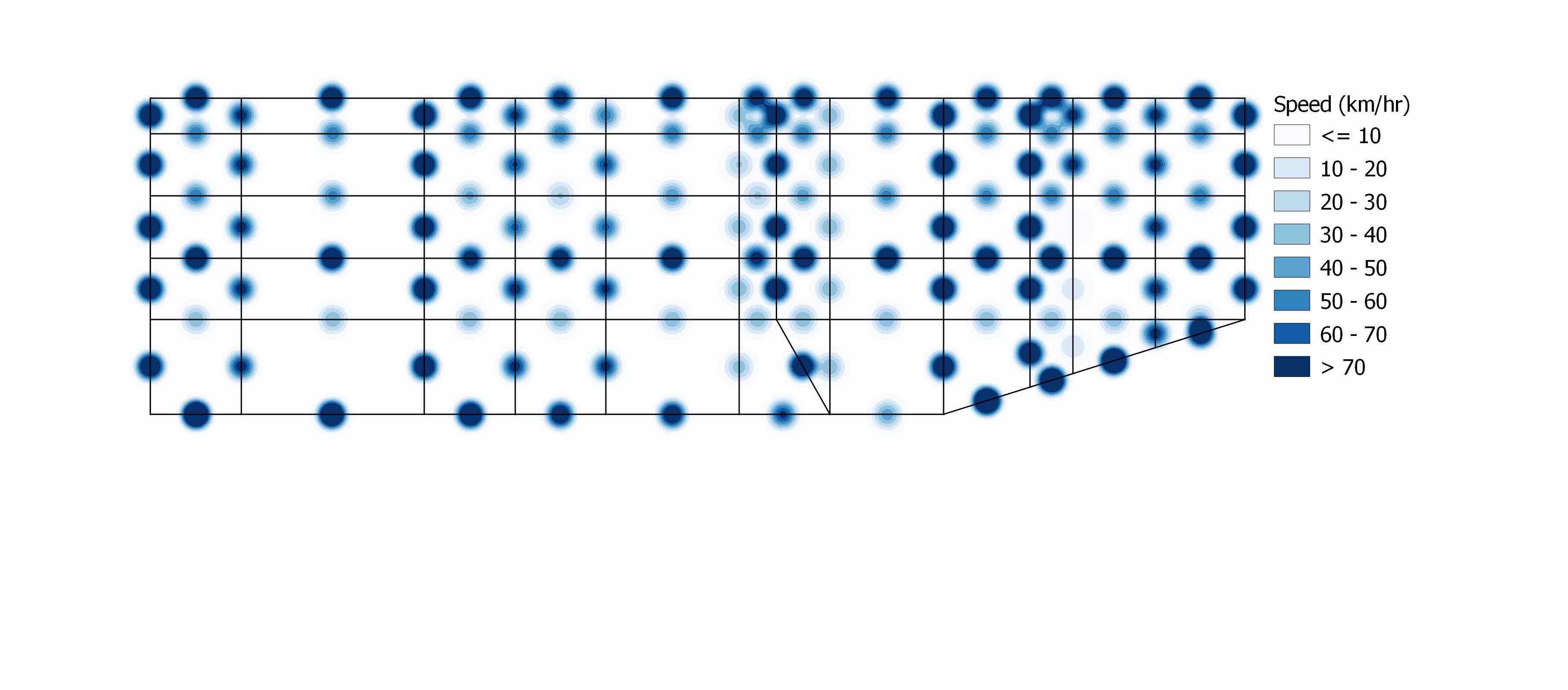} 
		\caption{Speed heatmap for S5}
		\label{speedheatmapTTstarplusGHG}
	\end{subfigure}%
	\begin{subfigure}{0.33\textwidth}
		\includegraphics[width=\textwidth]{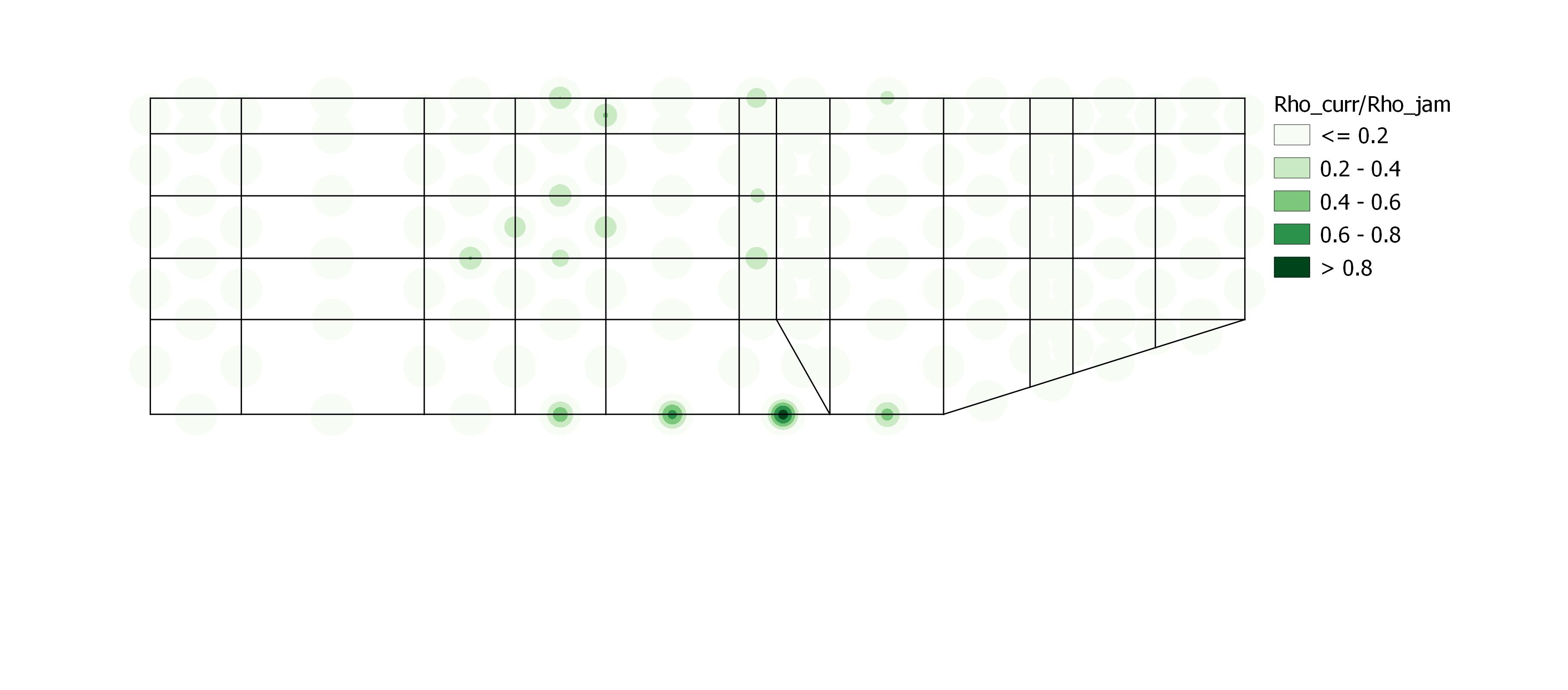}
		\caption{Density heatmap of S5}
		\label{DensityheatmapofTTstarplusGHG}
	\end{subfigure}%
	\begin{subfigure}{0.33\textwidth}
		\includegraphics[width=\textwidth]{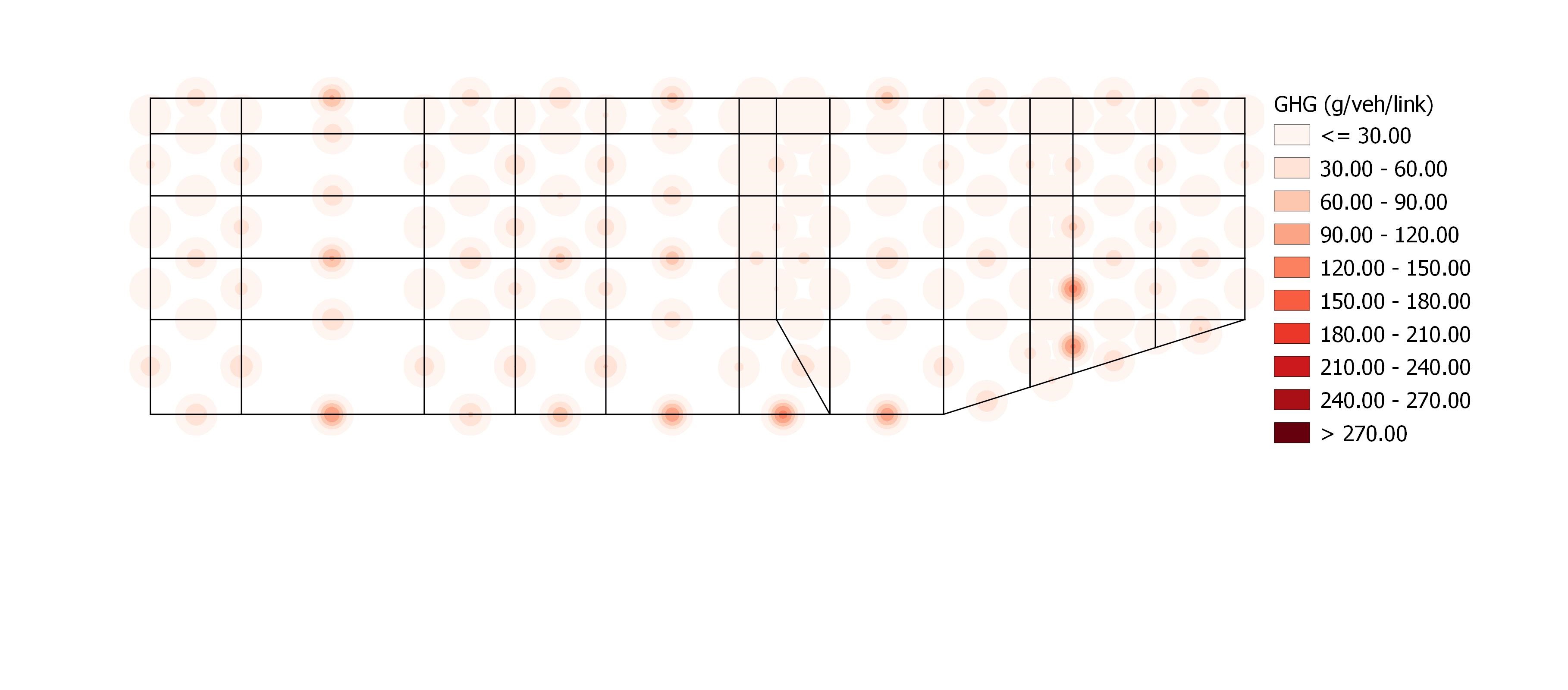}
		\caption{GHG heatmap of S5}
		\label{TTstarheatmapofGHG}
	\end{subfigure}
	\caption{Speed, density and GHG heatmaps of S5}
	\label{speed, density and GHG heatmap of TTstarplusGHG}
\end{figure}

\subsection{Air quality impact}
To further illustrate the impact of different routing algorithms on urban air quality, NO$_2$ concentrations are estimated through an atmospheric air pollution dispersion model SIRANE \cite{soulhac2011model}. SIRANE simulates the air pollution concentrations, taking link-based emission rate (gram/second) as the emission line source. Due to its consideration of urban geometry (such as street width and near-road building height), it is especially suitable for study domains where street canyon effect exists. Segment NO$_x$ emissions are estimated as the form of NO and NO$_2$ based on NO/NO$_2$ ratio in MOVES database; and the emission rate on one link as the line emission source in SIRANE during one evaluation interval equals total NO$_2$ emissions generated on that link in this period divided by total time duration. Detailed SIRANE parameter values adopted in this paper can be found in \cite{tu2019quantifying}. Since NO and NO$_2$ have similar spatial distributions under each scenario, in this paper, only NO$_2$ concentrations are presented. Results of scenarios S2--S5 are evaluated. Scenario S1 is omitted from this analysis because the aim here is to compare the effects of different E2ECAV based routing objectives on air quality. For each scenario, the whole simulation duration is extended to two hours and split into eight 15-minutes (900 seconds) intervals. It should be noted that the chemical reaction model adopted in SIRANE only includes a simplified photochemical reaction process. Therefore, the analysis in this section reveals only the impact from different routing scenarios, instead of a precise air pollution estimation from complex atmospheric reactions.  

\begin{figure}[!ht]
	\begin{subfigure}{0.5\textwidth}
		\includegraphics[width=\linewidth, height=4.5cm]{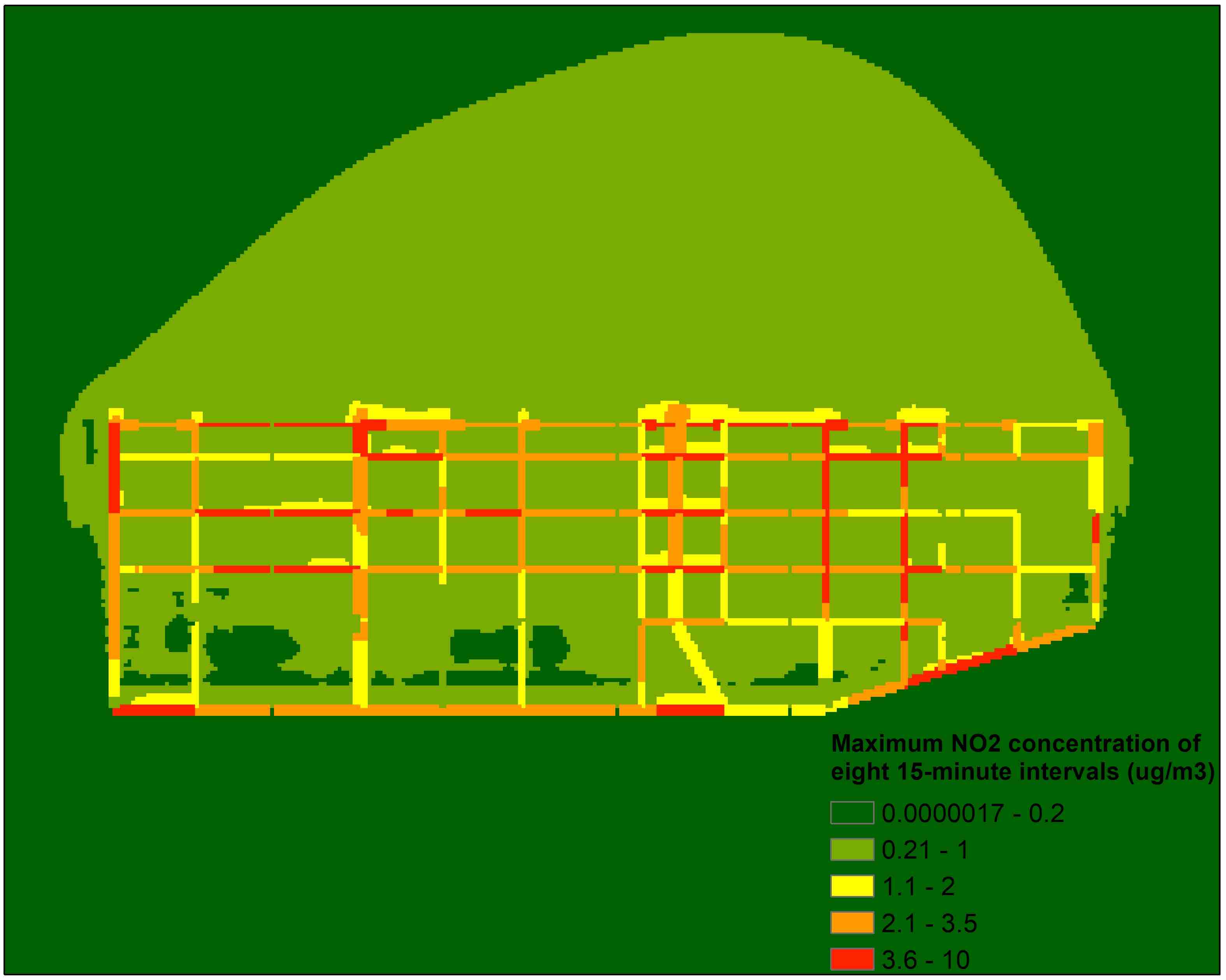} 
		\caption{S2}
		\label{Max S2 NO2 concentration}
	\end{subfigure}%
	\hfill
	\begin{subfigure}{0.5\textwidth}
		\includegraphics[width=\linewidth, height=4.5cm]{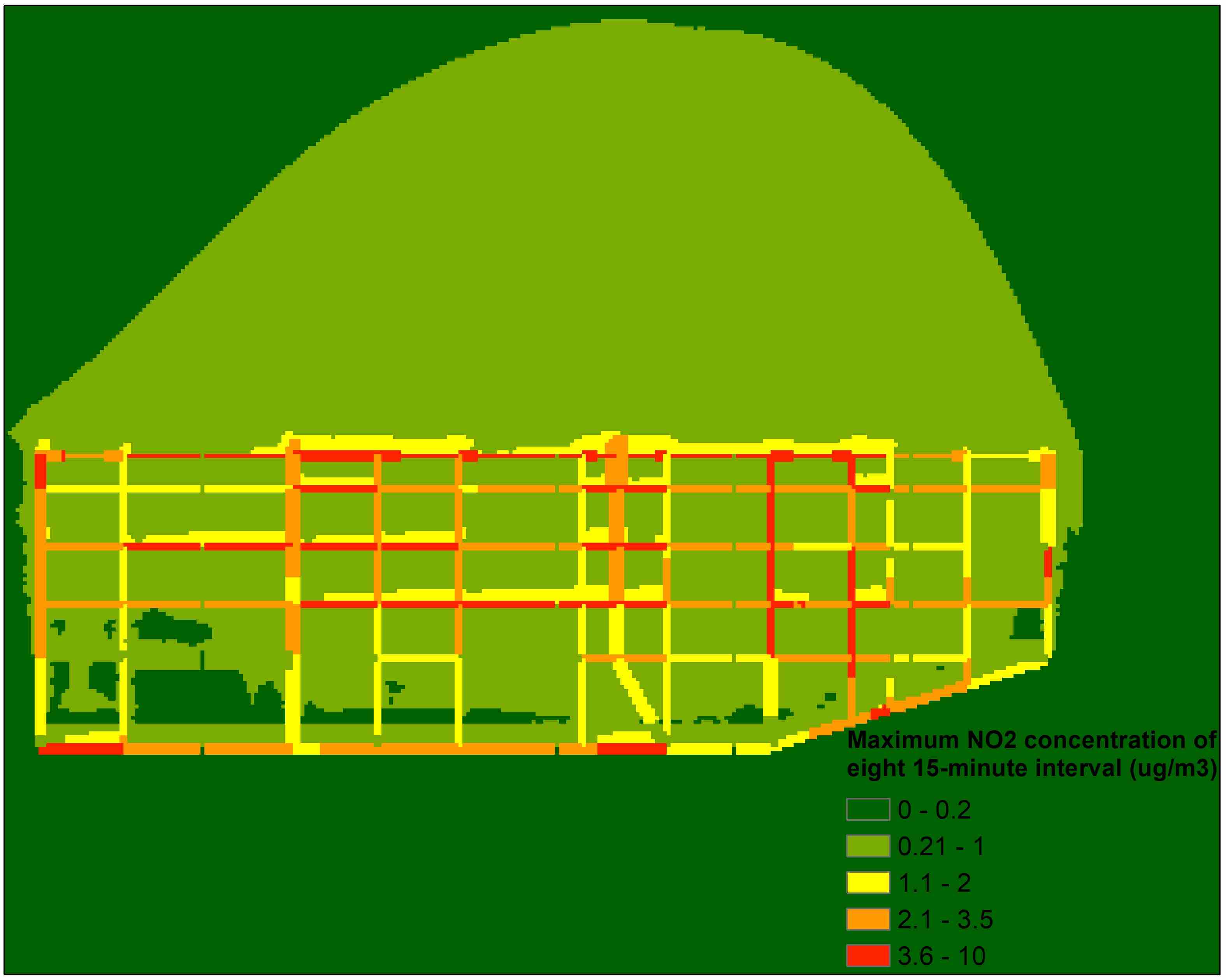} 
		\caption{S3}
		\label{Max S3 NO2 concentration}
	\end{subfigure}%
	\vskip\baselineskip
	\begin{subfigure}{0.5\textwidth}
		\includegraphics[width=\linewidth, height=4.5cm]{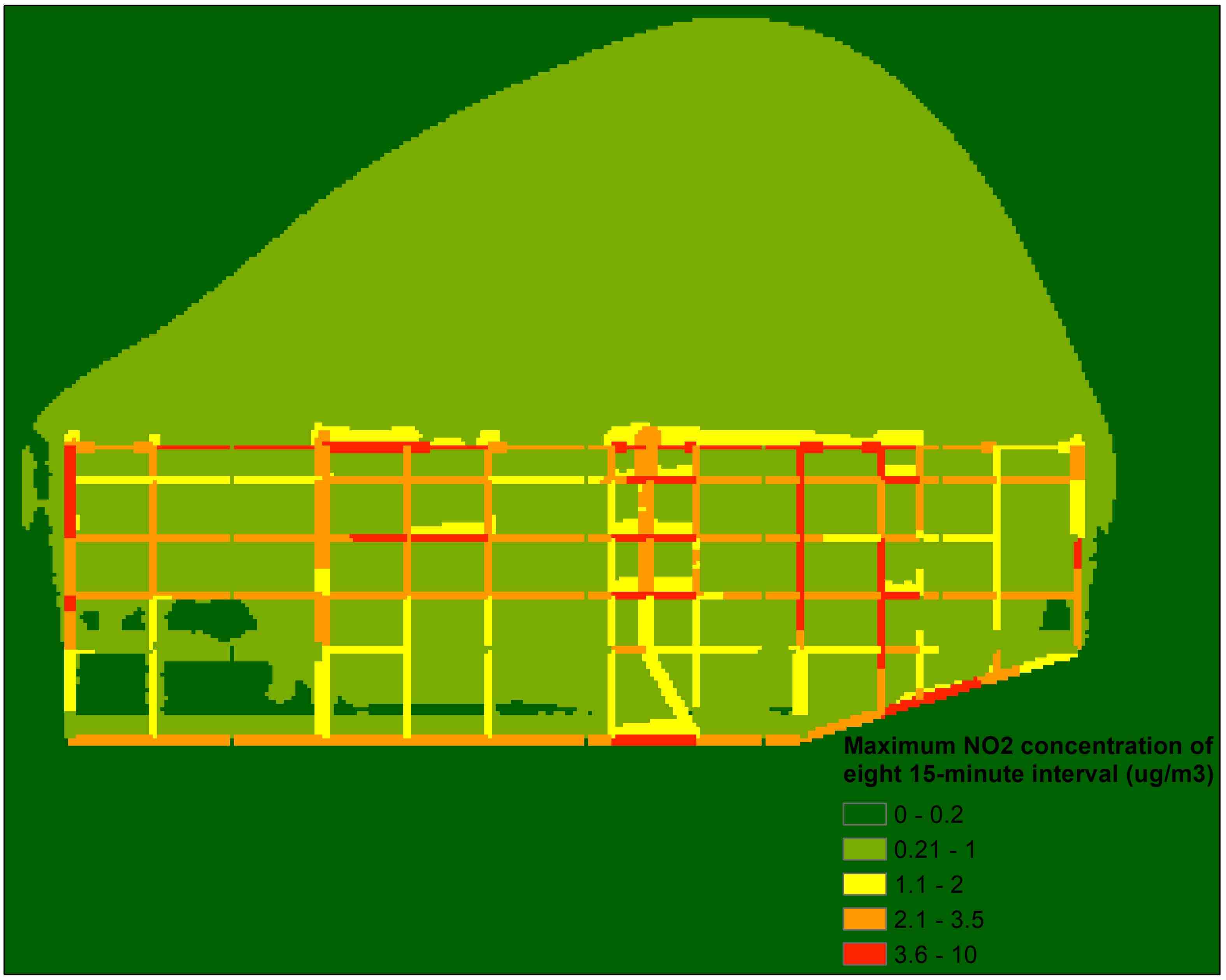}
		\caption{S4}
		\label{Max S4 NO2 concentration}
	\end{subfigure}
	\hfill
	\begin{subfigure}{0.5\textwidth}
		\includegraphics[width=\linewidth, height=4.5cm]{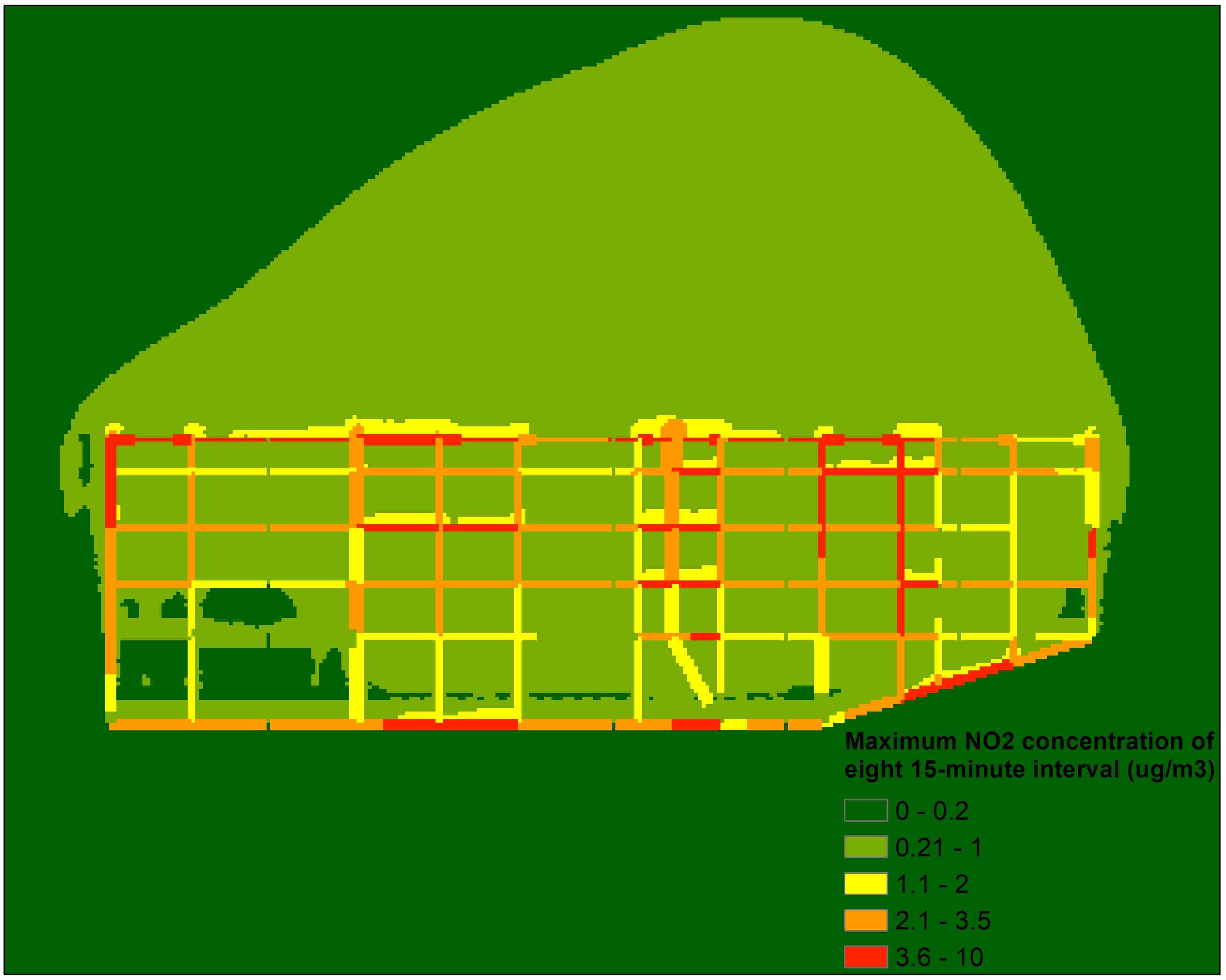}
		\caption{S5}
		\label{Max S5 NO2 concentration}
	\end{subfigure}
	\caption{2-hour NO$_2$ maximum concentrations of routing scenarios S2--S5}
	\label{2-hour NO$_2$ maximum concentrations of four routing scenarios}
\end{figure}

\begin{figure}[!ht]
	\begin{subfigure}{0.5\textwidth}
		\includegraphics[width=\linewidth, height=4.5cm]{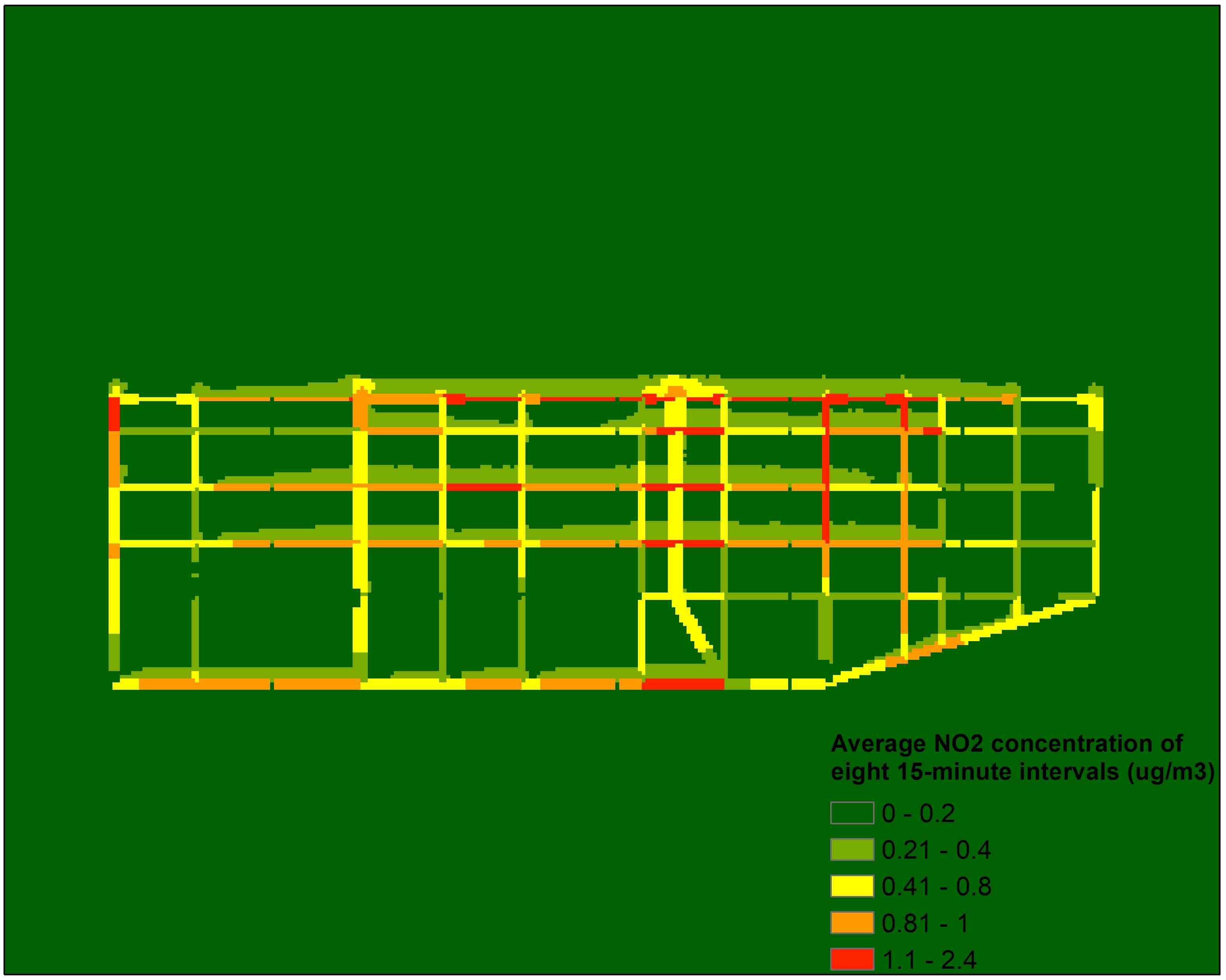} 
		\caption{S2}
		\label{avg S2 NO2 concentration}
	\end{subfigure}%
	\hfill
	\begin{subfigure}{0.5\textwidth}
		\includegraphics[width=\linewidth, height=4.5cm]{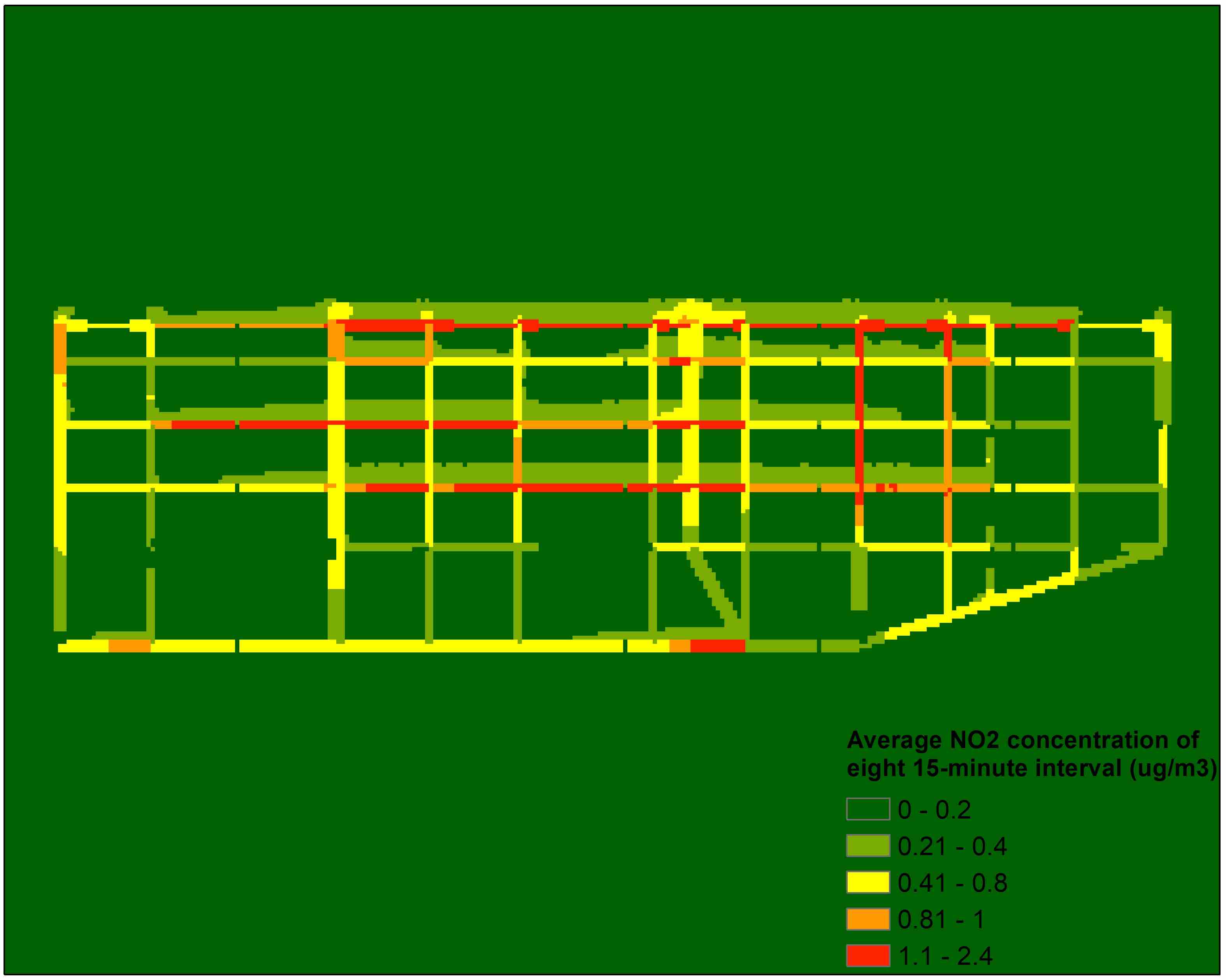} 
		\caption{S3}
		\label{avg S3 NO2 concentration}
	\end{subfigure}%
	\vskip\baselineskip
	\begin{subfigure}{0.5\textwidth}
		\includegraphics[width=\linewidth, height=4.5cm]{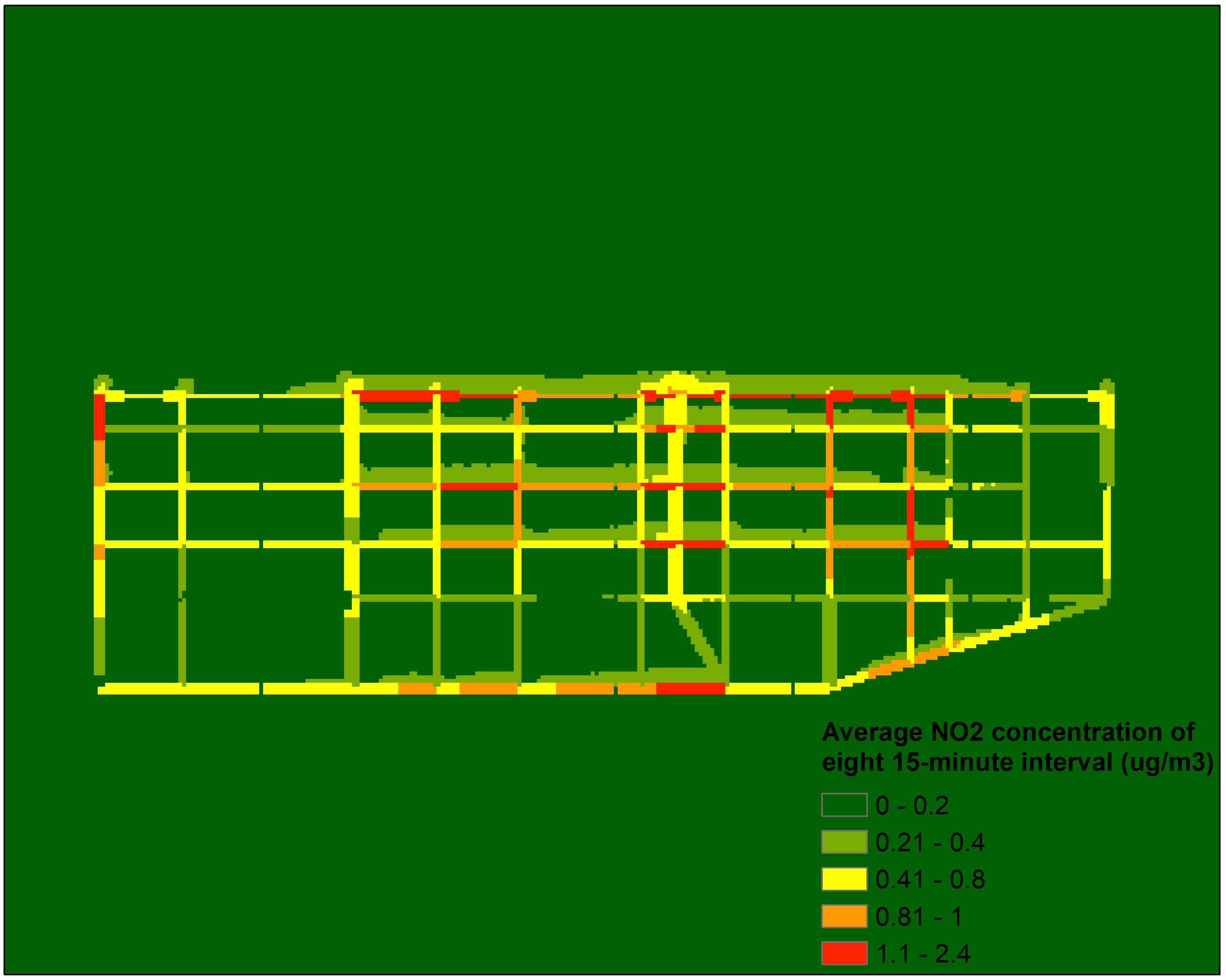}
		\caption{S4}
		\label{avg S4 NO2 concentration}
	\end{subfigure}
	\hfill
	\begin{subfigure}{0.5\textwidth}
		\includegraphics[width=\linewidth, height=4.5cm]{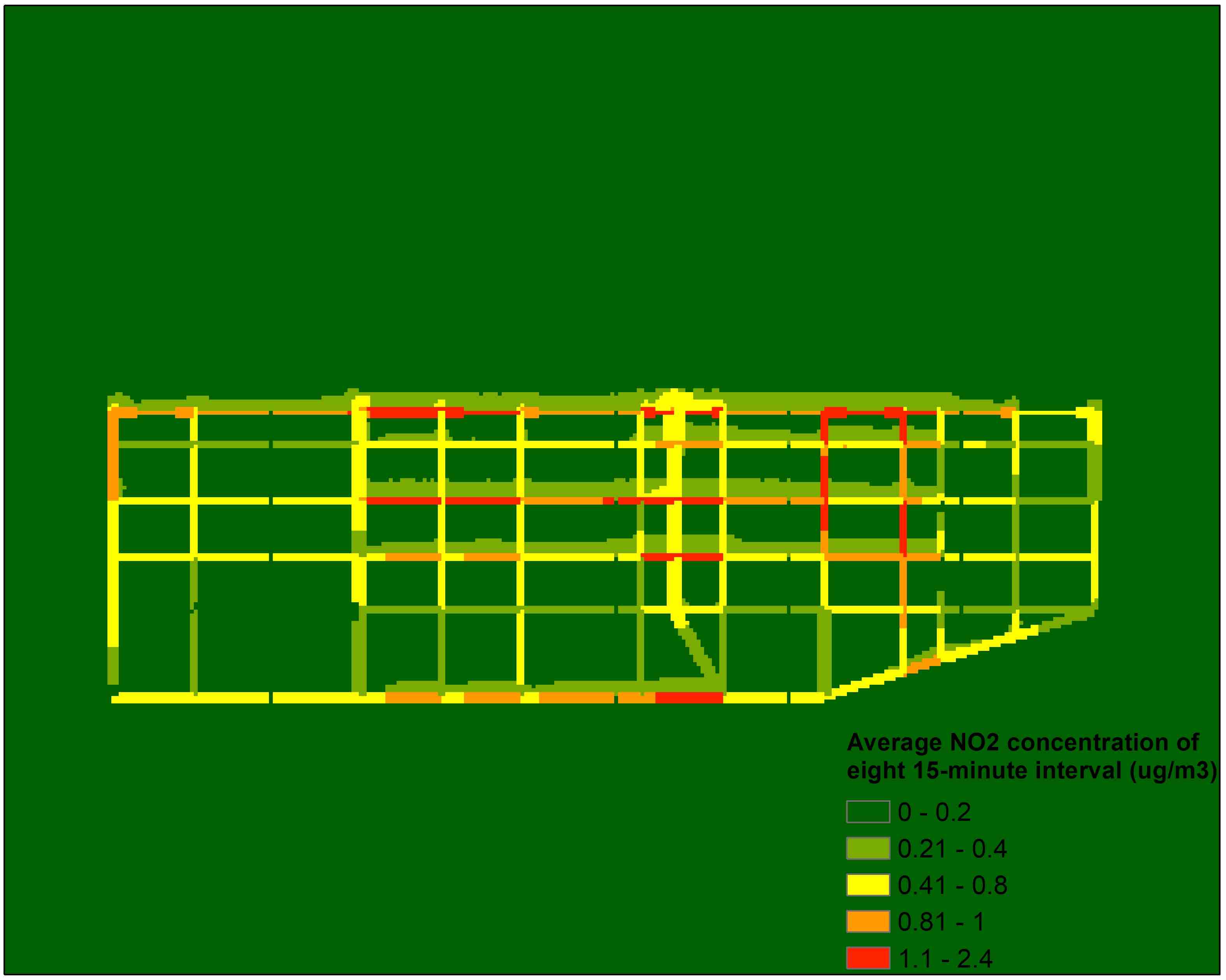}
		\caption{S5}
		\label{avg S5 NO2 concentration}
	\end{subfigure}
	\caption{2-hour NO$_2$ average concentrations of routing scenarios S2--S5}
	\label{2-hour NO$_2$ average concentrations of four routing scenarios}
\end{figure}

Figure \ref{2-hour NO$_2$ maximum concentrations of four routing scenarios} and Figure \ref{2-hour NO$_2$ average concentrations of four routing scenarios} demonstrate heatmaps for the maximum and average NO$_2$ concentrations among eight intervals in each 10mx10m grid, respectively. Including travel time in the objective function results in a similar spatial distribution of NO$_2$ concentration hot-spots, as it can be observed from S2, S4, and S5. In contrast, S3 where only GHG minimization is the objective, has a higher NO$_2$ concentration level on King St. and Adelaide St. Meanwhile, Front St. has lower NO$_2$ concentrations under S3. These results are aligned with the analysis on speed, density, and GHG emissions in Figure \ref{speed, density and GHG heatmap of GHG}. This indicates that only minimizing GHG emissions tends to re-route vehicles to streets with similar geometric design but lower speed limit. In contrast, both S4 and S5, with $TT$ based routing, are able to reduce NO$_2$ concentration level on King St. and Adelaide St. by re-routing the traffic to Front St. This leads to not only less emissions, but better utilization of the available infrastructure to reduce congestion. This also shows that although routing targeting only GHG emissions has the ability to reduce both GHG and NO$_x$ emissions, it may lead to inefficient utilization of the road network. 
Another finding from Figure \ref{2-hour NO$_2$ maximum concentrations of four routing scenarios} and Figure \ref{2-hour NO$_2$ average concentrations of four routing scenarios} is that the \textcolor{black}{idling penalty} might be a dominant factor that leads to a lower NO$_2$ concentration level. As we compare S2 with S4, the number of maximum NO$_2$ hot-spots decreases, and average NO$_2$ concentrations on the majority of links decrease from 0.8-1 $\mu$g in S2 to 0.41-0.8 $\mu$g in S4. However, S4 and S5 that consider \textcolor{black}{idling penalty}, have an almost identical spatial distribution of NO$_2$ hot-spots, with only minor differences in values.

\newpage

\newpage
\section{Conclusion}
\label{S:6}
Road transportation is one of the leading contributors to environmental pollution. This study develops a multi-objective eco-routing system using a real-time End-to-End CAV (E2ECAV) routing scheme \cite{FarooqBilalandDjavadian},  with the aim to simultaneously minimize travel time, GHG and NO$_x$ emissions. Particularly, two eco-routing strategies and associated objective functions are proposed. In total, five scenarios are developed based on these objective functions and vehicle routing schemes. These scenarios are applied to downtown Toronto network in an in-house agent-based traffic simulation platform. The performance measures used are: travel time, vehicle kilometre travelled, traffic GHG and NO$_x$ emissions. MOtor Vehicle Emission Simulator (MOVES) developed by the USEPA is employed to estimate the emission of a vehicle at a microscopic level. 

The results further emphasized the importance of taking into account the real-time estimates of emissions, travel time, delay, as well as the non-monotonic relationship between GHG and speed. The multi-objective eco-routing scenario that incorporated GHG emissions, travel time and delay in the objective function led to the highest reduction in the experienced travel time (40\%) and GHG emissions (43\%) when compared with the base case of human driven vehicles \textcolor{black}{with pre-trip dynamic routing}. By including the travel time in the objective function and having access to the real-time traffic conditions on the network, the routing strategy has the ability to periodically re-distribute traffic to the less congested parts of the network. Accounting for the average \textcolor{black}{idling penalty} at the downstream intersections, the strategy is able to prevent longer delays and dissipate the current queues faster. The inclusion of GHG emissions in the routing objective function resulted in further travel time reduction. Since the relationship between GHG and speed is quasi-convex, having GHG as an objective ensured a balance in distributing traffic on the network and better utilizing the available infrastructure by maintaining optimal speed on each link (not too low and not too high). 

In addition to travel time and GHG emission, the results showed that the proposed multi-objective routing could potentially reduce NO$_x$ emissions level by 18.5\%, despite the fact that it is not explicitly incorporated in the objective function. This significant improvement was achievable due the multi-objective eco-routing's ability to indirectly address the main factors affecting NO$_x$ (i.e. long travel time and high speed). By minimizing travel time and GHG it ensures that the vehicles are not staying too long on the network and also not travelling at high speeds. Furthermore, in comparison to the scenarios where either GHG emissions or travel time are optimized, incorporating \textcolor{black}{idling penalty} into the routing not only reduced the GHG emissions more effectively, but because of more efficient road capacity utilization, it also mitigated the near-road air quality along major corridors. The air quality improvements from the \textcolor{black}{idling penalty}-integrated routing algorithm would further benefit pedestrians on sidewalks.

In future research, the constrained optimization problem can be further refined. Applying weights to various competing terms in the multi-objective function would manage the differences in the scale of variables better. In this study we assumed 100\%CAVs on the network. Investigating the impact of different Market Penetration Rates (MPRs) of CAVs should be considered to illustrate the efficiency of the evaluated eco-routing strategies. In the routing scheme used here, traffic characteristics (travel time, \textcolor{black}{idling penalty} and emissions level) are updated every 60sec by averaging over the past interval. In the case of recurrent traffic conditions, the aforementioned updating method provides a good indicator of a prevailing condition on each link. However, in future studies more advanced prediction of the change of emissions and traffic parameters using machine learning methods can be brought in---especially in the case of non-recurrent traffic conditions. In addition, different algorithm settings can be explored in future research. For example, NO$_x$ emissions, which are more sensitive to aggressive driving behaviors like sharp acceleration and high speed, may lead to different results as GHG-based optimization. Therefore, the influence of NO$_x$-based optimization or a combined measure of different types of emissions can be evaluated. The negative impact of NO$_x$ pollutant on public health should be taken into consideration. Future studies will also look into minimization of the effect of near-road fine particles in urban areas. 

\textcolor{black}{It is worth mentioning that although the proposed multi-objective eco-routing was implemented on a distributed traffic management system, the findings are agnostic the routing scheme as long as it can exploit vehicular connectivity. In this study we proposed certain objective functions that need to be
	optimized and outlined the type of information that is needed in these objective functions. As
	long as a routing scheme is capable of providing the required information, optimizing the
	objective functions, and communicating with individual vehicles, our strategies will work.
	In the future work, we will explore the implementation of the proposed strategies on other routing schemes and communication schemes (e.g. V2V). }

\section*{Acknowledgement}
This research is funded by NSERC Canada Research Chair on Disruptive Transportation Technologies and Services, Ontario Early Researcher Award, and Ryerson University. We greatly appreciate Nikki Hera-Farooq's time and effort in helping us improve the quality of this manuscript.

\appendix




\section{Choice of routing for HDV}
\label{hdv_routing}
\begin{figure}[!h]
	\begin{center}
		\includegraphics[width=4in]{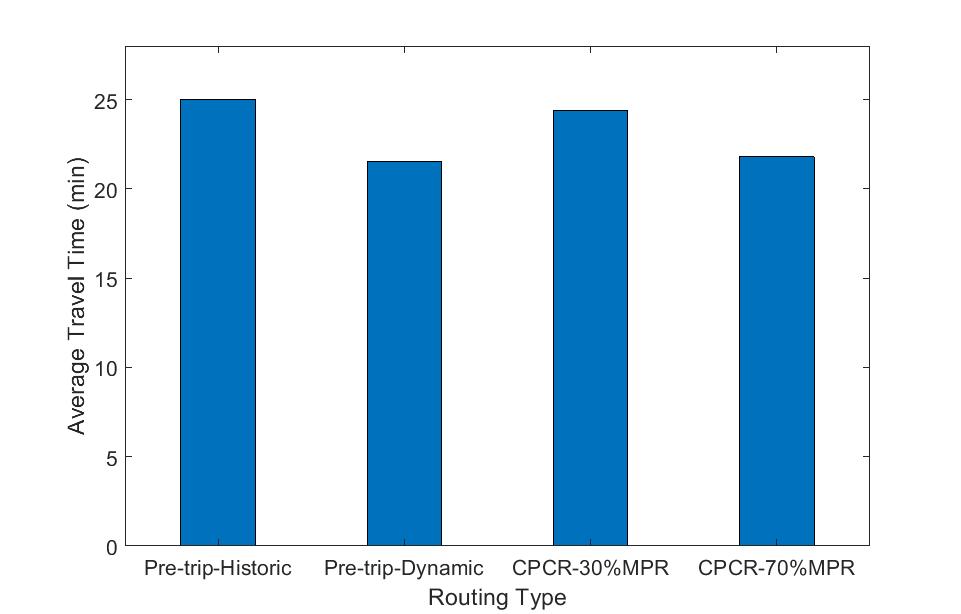}
		\caption{Travel time comparison for HDV routing}
		\label{hdvdynamic_vs_cpcr}
	\end{center}
\end{figure}

\textcolor{black}{
	In the case of scenario 1 ( HDV) a pre-trip dynamic shortest path is selected at the time of entry in the network based on current traffic conditions with no en-route routing.
	The reason for HDV pre-trip dynamic routing being selected as a based case is that we compared the performance of the above-mentioned HDV routing with \emph{Crowd-sourced Probe based Central Routing (CPCR)} (e.g. WAZE) under 30\% and 70\% market penetration and the results obtained showed that HDV based on pre-trip dynamic routing performed as well as crowd-sourcing service with 70\% market penetration rate as presented in Figure \ref{hdvdynamic_vs_cpcr}, based on this the dynamic pre-routing was selected as a benchmark. Figure \ref{hdvdynamic_vs_cpcr} presents the average travel times for HDVs under three different routing schemes, namely:  pre-trip historic routing, pre-trip dynamic routing and CPCR routing. The reason HDV dynamic pre-trip routing performed as well as CPCR-70\% and outperformed CPCR-30\% is that although in the case of HDV with dynamic pre-trip routing,  en-route routing is not available the traffic information used to develop the pre-trip routing is updated more frequently and is more reliable. }

\newpage
\bibliographystyle{abbrv}
\bibliography{bibliography.bib}

\begin{thebibliography}{10}

\bibitem{AlMallah2017a}
R.~{Al Mallah}, A.~Quintero, and B.~Farooq.
\newblock {Distributed Classification of Urban Congestion Using VANET}.
\newblock {\em IEEE Transactions on Intelligent Transportation Systems},
  18(9):2435--2442, 2017.

\bibitem{AlMallah2017b}
R.~{Al Mallah}, A.~Quintero, and B.~Farooq.
\newblock {Cooperative Evaluation of the Cause of Urban Traffic Congestion via
  VANET}.
\newblock Technical report, Polytechnique Montreal and Ryerson University,
  2018.

\bibitem{alfaseeh2018impact}
L.~Alfaseeh, S.~Djavadian, and B.~Farooq.
\newblock Impact of distributed routing of intelligent vehicles on urban
  traffic.
\newblock In {\em 2018 IEEE International Smart Cities Conference (ISC2)},
  pages 1--7. IEEE, 2018.

\bibitem{alfaseeh2019eco}
L.~Alfaseeh, S.~Djavadian, and B.~Farooq.
\newblock Multi-objective eco-routing in a distributed routing framework.
\newblock In {\em 2019 IEEE International Smart Cities Conference (ISC2)},
  pages 1--7. IEEE, 2019.

\bibitem{alfaseeh2019multi}
L.~Alfaseeh and B.~Farooq.
\newblock Multi-factor taxonomy of eco-routing models and future outlook.
\newblock {\em Journal of Sensors}, pages 1--18, 2020.

\bibitem{Anagnostopoulou2018}
E.~Anagnostopoulou, E.~Bothos, B.~Magoutas, J.~Schrammel, and G.~Mentzas.
\newblock {Persuasive technologies for sustainable mobility: State of the art
  and emerging trends}.
\newblock {\em Sustainability (Switzerland)}, 2018.

\bibitem{andersen2013ecotour}
O.~Andersen, C.~S. Jensen, K.~Torp, and B.~Yang.
\newblock Ecotour: Reducing the environmental footprint of vehicles using
  eco-routes.
\newblock In {\em Mobile Data Management (MDM), 2013 IEEE 14th International
  Conference on}, volume~1, pages 338--340. IEEE, 2013.

\bibitem{aziz2012integration}
H.~A. Aziz and S.~V. Ukkusuri.
\newblock Integration of environmental objectives in a system optimal dynamic
  traffic assignment model.
\newblock {\em Computer-Aided Civil and Infrastructure Engineering},
  27(7):494--511, 2012.

\bibitem{GHG_cost}
T.~W. Bank, 2019.
\newblock [Online; accessed 10-April-2019].

\bibitem{bektacs2011pollution}
T.~Bekta{\c{s}} and G.~Laporte.
\newblock The pollution-routing problem.
\newblock {\em Transportation Research Part B: Methodological},
  45(8):1232--1250, 2011.

\bibitem{Canada2016}
N.~R. Canada.
\newblock {Emission impacts resulting from vehicle idling}, 2016.

\bibitem{djavadian2018distributed}
S.~Djavadian and B.~Farooq.
\newblock Distributed dynamic routing using network of intelligent
  intersections.
\newblock In {\em ITS Canada Annual General Meeting Conference, Niagara Falls},
  2018.

\bibitem{DjavadianS.FarooqB.VasquezR.Yip2019}
S.~Djavadian, B.~Farooq, R.~Vasquez, and G.~Yip.
\newblock {Virtual Immersive Reality based Analysis of Behavioral Responses in
  Connected and Autonomous Vehicle Environment}.
\newblock In A.~{Goulias, K. Davis}, editor, {\em Mapping the Travel Behavior
  Genome}, chapter~27. Elsevier, 2019.

\bibitem{FarooqBilalandDjavadian}
B.~Farooq and S.~Djavadian.
\newblock {Distributed Traffic Management System with Dynamic End-to-End
  Routing, U.S. Provisional Pat. Ser. No. 62/865,725}, 2019.

\bibitem{grote2016including}
M.~Grote, I.~Williams, J.~Preston, and S.~Kemp.
\newblock Including congestion effects in urban road traffic co2 emissions
  modelling: do local government authorities have the right options?
\newblock {\em Transportation Research Part D: Transport and Environment},
  43:95--106, 2016.

\bibitem{guo2013evaluation}
L.~Guo, S.~Huang, and A.~W. Sadek.
\newblock An evaluation of environmental benefits of time-dependent green
  routing in the greater buffalo--niagara region.
\newblock {\em Journal of Intelligent Transportation Systems}, 17(1):18--30,
  2013.

\bibitem{Kattan2012}
L.~Kattan, M.~Mousavi, B.~Far, C.~Harschnitz, A.~Radmanesh, and S.~Saidi.
\newblock {Microsimulation Evaluation of the Potential Impacts of
  Vehicle-to-Vehicle Communication (V2V) in Disseminating Warning Information
  under High Incident Occurrence Conditions}.
\newblock {\em International Journal of Intelligent Transportation Systems
  Research}, 10(3):137--147, 2012.

\bibitem{Liu2019}
H.~Liu, R.~Guensler, H.~Lu, Y.~Xu, X.~Xu, and M.~O. Rodgers.
\newblock {MOVES-Matrix for high-performance on-road energy and running
  emission rate modeling applications}.
\newblock {\em Journal of the Air and Waste Management Association}, 2019.

\bibitem{long2016link}
J.~Long, J.~Chen, W.~Szeto, and Q.~Shi.
\newblock Link-based system optimum dynamic traffic assignment problems with
  environmental objectives.
\newblock {\em Transportation Research Part D: Transport and Environment},
  60:56--75, 2016.

\bibitem{Meshkani2019}
S.~M. Meshkani, S.~Djavadian, and B.~Farooq.
\newblock {A Decentralizd Shared CAV system Design and Application}.
\newblock In {\em the Tenth Triennial Symposium on Transportation Analysis
  (TRISTAN)}, Hamilton Island, Australia, 2019.

\bibitem{Metrolinx2008}
Metrolinx.
\newblock {Costs of Roads Congestion in the Greater Toronto and Hamilton Area:
  Impact and Cost Benefit Analysis of the Metrolinx Draft Regional
  Transportation Plan}.
\newblock Technical report, 2008.

\bibitem{Michon1985}
J.~A. Michon.
\newblock {A Critical View of Driver Behavior Models: What Do We Know, What
  Should We Do?}
\newblock In L.~{L. Evans} and R.~Schwing, editors, {\em Human Behavior and
  Traffic Safety}, pages 485--524. New York: Plenum Press, 1985.

\bibitem{MTO2016}
{Ministry of Transportation Ontario}.
\newblock Statistics and management reporting, 2016.

\bibitem{Noori2013}
H.~Noori and M.~Valkama.
\newblock {Impact of VANET-based V2X communication using IEEE 802.11p on
  reducing vehicles traveling time in realistic large scale urban area}.
\newblock In {\em 2013 International Conference on Connected Vehicles and Expo,
  ICCVE 2013 - Proceedings}, pages 654--661, 2013.

\bibitem{patil2016emission}
G.~R. Patil.
\newblock Emission-based static traffic assignment models.
\newblock {\em Environmental Modeling \& Assessment}, 21(5):629--642, 2016.

\bibitem{Sawyer2019}
P.~Sawyer.
\newblock {How Waze is using data pacts, beacons, and carpools to win over
  cities}.
\newblock {\em Venturebeat}, 2019.

\bibitem{soulhac2011model}
L.~Soulhac, P.~Salizzoni, F.-X. Cierco, and R.~Perkins.
\newblock The model sirane for atmospheric urban pollutant dispersion; part i,
  presentation of the model.
\newblock {\em Atmospheric environment}, 45(39):7379--7395, 2011.

\bibitem{sun2015stochastic}
J.~Sun and H.~X. Liu.
\newblock Stochastic eco-routing in a signalized traffic network.
\newblock {\em Transportation Research Procedia}, 7:110--128, 2015.

\bibitem{Transport&Environment}
{Transport {\&} Environment}.
\newblock {Don't Breath Here: beware the invisible killer}.
\newblock Technical report, 2015.

\bibitem{Treiber2000}
M.~Treiber, A.~Hennecke, and D.~Helbing.
\newblock {Congested traffic states in empirical observations and microscopic
  simulations}.
\newblock {\em Physical Review E - Statistical Physics, Plasmas, Fluids, and
  Related Interdisciplinary Topics}, 62(2):1805--1824, 2000.

\bibitem{tu2019quantifying}
R.~Tu, L.~Alfaseeh, S.~Djavadian, B.~Farooq, and M.~Hatzopoulou.
\newblock Quantifying the impacts of dynamic control in connected and automated
  vehicles on greenhouse gas emissions and urban no2 concentrations.
\newblock {\em Transportation Research Part D: Transport and Environment},
  73:142--151, 2019.

\bibitem{UnitedStatesEnvironmentalProtectionAgency2019}
{United States Environmental Protection Agency}.
\newblock {MOVES2014b: Latest Version of MOtor Vehicle Emission Simulator
  (MOVES)}, 2019.

\bibitem{usdot2014}
USDOT.
\newblock The value of travel time savings: Departmental guidance for
  conducting economic evaluations.
\newblock Technical report, US. Department of Transportation, 2014.

\bibitem{yang2006modeling}
X.~Yang and W.~W. Recker.
\newblock Modeling dynamic vehicle navigation in a self-organizing,
  peer-to-peer, distributed traffic information system.
\newblock {\em Journal of intelligent transportation Systems}, 10(4):185--204,
  2006.

\end{thebibliography}
\end{document}